\documentclass{article}
\newif\ifclean
\cleanfalse
\cleantrue
\usepackage{arxiv}
\usepackage[letterpaper]{geometry}
\usepackage{graphicx}
\usepackage{subcaption}
\usepackage{amsmath,amssymb}
\usepackage{algorithm}
\usepackage{algpseudocode}
\graphicspath{{{./}}}

\ifclean
\else
\usepackage{soul}
\usepackage{xcolor}
\newcommand{\red}[1]{\textcolor{red}{{#1}}}

\fi

\newcommand{\eref}[1]{{Eq.\,\eqref{#1}}}
\newcommand{\tref}[1]{{Tab.\,\ref{#1}}}
\newcommand{\fref}[1]{{Fig.\,\ref{#1}}}
\newcommand{\cref}[1]{{Ref.\,\cite{#1}}}
\newcommand{\sref}[1]{{Sec.\,\ref{#1}}}

\newcommand{\parameters}{{\boldsymbol{\theta}}}
\newcommand{\outpt}{{\mathsf{y}}}
\newcommand{\modeloutpt}{{\tilde{\mathsf{y}}}}

\newcommand{\inputs}{{\mathsf{X}}}
\newcommand{\inpt}{{\mathsf{x}}}
\newcommand{\outputs}{{\mathsf{Y}}}

\newcommand{\residuals}{{\mathsf{R}}}
\newcommand{\sensitivities}{{\mathsf{A}}}
\newcommand{\hessian}{{\mathsf{H}}}
\newcommand{\model}{{\mathbf{m}}}
\newcommand{\data}{{\mathbf{d}}}
\newcommand{\dataset}{{\mathsf{D}}}
\newcommand{\noise}{{\boldsymbol{\eta}}}
\newcommand{\pvariance}{{\boldsymbol{\xi}}}
\newcommand{\varparameters}{{\boldsymbol{\vartheta}}}

\newcommand{\prob}{{p}}
\newcommand{\likelihood}{{L}_{\mathsf{D}}}
\newcommand{\entropy}{{H}}

\newcommand{\EXP}{{\mathbb{E}}}
\newcommand{\VAR}{{\mathbb{V}}}
\newcommand{\KLD}{\operatorname{KLdiv}}
\newcommand{\CDF}{\operatorname{CDF}}
\newcommand{\Sigmab}{{\boldsymbol{\Sigma}}}
\newcommand{\Thetab}{{\boldsymbol{\Theta}}}
\newcommand{\thetab}{{\boldsymbol{\theta}}}
\newcommand{\varthetab}{{\boldsymbol{\vartheta}}}
\newcommand{\phib}{{\boldsymbol{\phi}}}
\newcommand{\partialb}{{\boldsymbol{\partial}}}

\newcommand{\objective}{{\Phi}}
\newcommand{\argmin}{\operatorname{argmin}}
\newcommand{\Argmin}[1]{\underset{{#1}}{\operatorname{argmin}}}

\newcommand{\As}{\mathsf{A}}
\newcommand{\Bs}{\mathsf{B}}

\newcommand{\Ds}{\mathsf{D}}
\newcommand{\Hs}{\mathsf{H}}

\newcommand{\Is}{\mathsf{I}}
\newcommand{\Ls}{\mathsf{L}}
\newcommand{\Ns}{\mathsf{N}}

\newcommand{\Us}{\mathsf{U}}

\newcommand{\Ws}{\mathsf{W}}

\newcommand{\Xs}{\mathsf{X}}
\newcommand{\Ys}{\mathsf{Y}}

\newcommand{\Rs}{\mathsf{R}}
\newcommand{\fb}{\mathbf{f}}

\newcommand{\zerob}{\mathbf{0}}

\newcommand{\Ib}{\mathbf{I}}

\newcommand{\varSigmab}{\boldsymbol{\Sigma}}

\newcommand{\mub}{\boldsymbol{\mu}}

\newcommand{\zetab}{\boldsymbol{\zeta}}

\newcommand{\Nc}{\mathcal{N}}
\newcommand{\Wc}{\mathcal{W}}

\newcommand{\tr}{\operatorname{tr}}

\newcommand{\arctanh}{\operatorname{arctanh}}

\newcommand{\stress}{\boldsymbol{\sigma}}
\newcommand{\strain}{\boldsymbol{\epsilon}}
\newcommand{\yield}{\Upsilon}
\newcommand{\internalstate}{\boldsymbol{\kappa}}

\title{\bf Laplace approximation of the parameter posterior with embedded error for efficient uncertainty quantification of expensive model calibrations}
\renewcommand{\shorttitle}{\bf Laplace approximation of the posterior with embedded error}
\author{\textrm{R.E.Jones, K.N. Karlson, M.W. Kury, M. Khalil, and R.S. Tuminaro} \\[0.05in]
\textit{Sandia National Laboratories, Livermore, CA 94551}}
\ifclean
\else
\date{\red{\textbf{Draft} \today: do not distribute} }
\fi

\begin{document}
\maketitle

\begin{abstract}
In practical engineering, there is a clear need for uncertainty quantification (UQ) of models used in design because forward UQ puts probabilities on performance predictions.
However, without surrogate building, most forward engineering models are too expensive for traditional UQ methods like Markov chain Monte Carlo (MCMC) sampling.
Our goal is to provide inexpensive, approximate UQ that can be applied after an expensive traditional least squares fit, where the UQ handles both uncorrelated noise and parameter variance as in the \emph{embedded-error} treatment developed by Sargsyan, Najm, and Ghanem (2015).
After developing the basic theory to connect a Laplace-like approximation to Bayesian calibration, we compare the proposed method to \textit{black-box} variational inference in terms of cost and fidelity.
Then we demonstrate the method on real-world calibration of solid material models.
\end{abstract}


\raggedbottom

\section{Introduction}

Uncertainty quantification (UQ) provides confidence in simulations used in design and performance estimation.
It can also uncover limitations of the calibration process, such as insufficient data and overly correlated parameters.
However, the sampling-based methods, such as Markov chain Monte Carlo \cite{ghanem2017handbook}, are typically too computationally expensive (and temperamental in terms of convergence and hyperparameter settings) to deploy for complex material models embedded in the full simulations of experiments used for calibration.
Inexact methods, such as variational inference \cite{blei2017variational}, are more efficient but also can be prohibitively expensive for these scenarios, depending on the optimization technique \cite{bottou2012stochastic}.

\emph{Laplace approximation} is a form of variational inference in the sense that a surrogate posterior is formed via an optimization procedure.
It leads to Gaussian approximation of the observed outputs at the maximum a posteriori (MAP) point in parameter space.
Since it only requires a Hessian at the optimum found by deterministic calibration, it can be extremely model sample/evaluation efficient, especially for a small number of parameters.
It is a classical technique going back to Laplace's method to approximate integrals \cite{tierney1986accurate, mackay1992bayesian}.
Here, we fuse this approach with the \emph{embedded error} formulation \cite{sargsyan2015statistical} to accommodate both uncorrelated external noise and correlated internal variations that are evident in the experimental data \cite{romero2019propagating}.
The embedded error formulation allows the model parameters to be random variables, as opposed to the classical Kennedy-O'Hagan \cite{kennedy2001bayesian} external/additive error, which enables the formulation to capture correlated deviations from a mean response.
This allows the separation of external, uncorrelated errors from internal, correlated errors.
The first source of uncertainty is usually due to measurement resolution/precision limits while the second is an outcome of specimen (or environmental) variability.
When there is significant material variability, the measurement noise is typically of smaller magnitude with modern testing equipment.

The inexact nature of the Laplace approximation is not a hindrance in our applications since the main  goal is to cover likely responses in prediction and hence have credible estimates for performance evaluation.
Although not a \emph{goal oriented} \cite{duraisamy2012goal,wahal2019bimc} technique, it has some commonalities with this more sophisticated class of methods.

In the next section we put the proposed method in the larger context of uncertainty quantification.
Then, in \sref{sec:method}, we develop the theoretical and methodological basis for the proposed method while emphasizing the contribution of extending the Laplace method to embedded error and a straightforward linear algebra solution for the parameter covariance.
In \sref{sec:results} we explore the algorithmic properties of the method with simple examples before providing results for real-world material calibrations with diverse experimental data.
Finally, we conclude with a summary and avenues for future work in \sref{sec:conclusion}.

\section{Related work}
Uncertainty quantification (UQ) has a long history of development \cite{ghanem2017handbook}.
Although still the reference standard, \emph{Markov chain Monte Carlo} (MCMC) sampling has a slow $1/\sqrt{n}$ convergence, where $n$ is the number of samples.
Newer techniques, such as \emph{Hamiltonian Monte Carlo} (HMC) \cite{neal2011mcmc} and Langevin sampling \cite{roberts1996exponential,roberts1998optimal,welling2011bayesian}, which employ gradients of the score and a connection to statistical physics, can be more sample efficient but still are prohibitive for expensive model evaluations, such as those from design-level finite element simulations.
The cost of these methods has motivated optimization-based techniques, such as \emph{variational inference} (VI) \cite{blei2017variational} where a surrogate posterior is determined by optimization.
However, the more scalable \emph{mean field} VI, which uses a diagonal covariance structure, is well known to under-predict the extent of posterior probabilities \cite{wang2005inadequacy}.
Hybrid methods, like \emph{Stein VI} \cite{liu2016stein} and other ensemble/\emph{particle} methods \cite{liu2019understanding}, try to obtain the best qualities of sampling and optimization techniques.

In all these approaches the primary concern, beyond the forward simulation cost, is size of the parameter space, whereby a large number of parameters to explore leads to the \emph{curse of dimensionality} \cite{bellman1957dynamic} due to the exponential rise in the number of samples needed to characterize a parameter or output covariance structure.
This fundamental issue emphasizes the need for efficient UQ for highly parameterized, expensive models.
The formulation of an appropriate likelihood tying observed data to model predictions is also a challenge.
Kennedy and O'Hagan \cite{kennedy2001bayesian} proposed an additive, independent error model-observation discrepancy that has been widely adopted and is particularly well suited for representing (external) measurement noise.
The Kennedy and O'Hagan approach leads to a straightforward likelihood formulation.
In contrast,  Sargsyan et al. \cite{sargsyan2015statistical,sargsyan2019embedded} proposed a correlated \emph{embedded} error model whereby parameters are associated with random variables which induce a covariance structure on the \emph{pushed forward} model responses.%
\footnote{A push-forward response, in this context, is taking samples of the uncertain parameters and evaluating the model with them in order to construct the corresponding  ensemble of outputs.}
This approach suggests a number of possible likelihood formulations.
These sources of discrepancies can be combined, but this increases the complexity of the likelihood and the number of hyperparameters that need to be selected/determined.

\section{Calibration} \label{sec:method}

In this section, we provide the theory necessary to understand the foundations of the proposed method.
The classical Laplace approximation provides an analytical approximation to posterior probability using a quadratic approximation of the log posterior distribution.
This Gaussian surrogate is typically centered at the maximum a posteriori (point) estimate of the parameters (which corresponds to the maximum likelihood fit for an uninformative prior).
The precision matrix is given by the Hessian at the center, which is  typically further approximated by a Gauss-Newton outer product of the model sensitivities (which is related to the Fisher information).

In the following developments, we start with least squares calibration (of an expensive model) given data, connect this to Bayesian calibration and, finally, to Laplace approximation.
The use of an embedded error model is emphasized.

\paragraph{Data}
We consider experimental data that consists of input-output pairs $\dataset = (\inputs,\outputs)$, which are subject to external measurement noise.
In many applications, such as material characterization via mechanical testing, measurement noise is high frequency and uncorrelated with the underlying response.
On the other hand, the underlying response can have intrinsic/internal variation that systematically affects the data trends.

\paragraph{Least squares calibration}
In general, a model $\model$ has parameters $\parameters$ and produces output $\modeloutpt$ as a function of an input $\inpt \in \inputs$:
\begin{equation} \label{eq:model}
\modeloutpt = \model(\inpt ; \parameters) \ .
\end{equation}
Given measured (experimental) data $\dataset=(\inputs,\outputs)$, the model can be calibrated using the usual least squares optimization
\begin{equation} \label{eq:least_squares}
\parameters^*
= \argmin_\parameters \Phi(\parameters ; \dataset, \model)
\equiv \argmin_\parameters \frac{1}{2} \| \underbrace{ \outputs -\model(\inputs; \parameters) }_{\residuals} \|^2_{\Ws}
=  \argmin_\parameters \frac{1}{2}\residuals^T \Ws \residuals
\ ,
\end{equation}
based on the (weighted) sum squares error (SSE) $\Phi$ to obtain the best fit parameters $\parameters^*$.
Here, the objective $\Phi$ uses the residuals:
\begin{equation} \label{eq:residuals}
\residuals = \outputs - \model(\inputs; \parameters) \ ,
\end{equation}
weighted by $\Ws$ which quantifies the uncertainty in $\outputs$ by identifying it with the inverse covariance (\emph{precision}) matrix of the responses/outputs $\Ws = \Sigmab_\outputs^{-1}$, as in standard generalized least squares calibration.
Both $\Ws$ and $\Sigmab_\outputs$ are of dimensions $n_y \times n_y$ where $n_y$ is the length of $\outputs$; hence experimental repeats are handled by making both $\outputs$ and $\residuals$ two dimensional with dimensions $n_y \times n_r$ where $n_r$ is the number of (data) repeats.
The inputs $\Xs$ are the same dimensions as $\Ys$ but identical/duplicated in the repeat dimension.

\paragraph{Bayesian calibration}
The observed data $\dataset=(\inputs,\outputs)$ can be related to the model by adding (probabilistic) mechanisms for discrepancies
\begin{equation}  \label{eq:embedded_model}
\outpt \approx \model(\inpt ; \bar{\parameters} + \pvariance) + \noise \quad \text{for} \ \ \inpt \in \inputs \ ,
\end{equation}
where we have output data $\outpt \in \outputs$ that displays (relatively small) uncorrelated measurement error (\emph{noise}) $\noise$ and correlated model inconsistencies (\emph{discrepancy}) which we associate with $\pvariance$, as mentioned.
Both $\noise$ and $\pvariance$ are considered (independent) random variables.
We associate $\noise$ with uncorrelated, irreducible measurement noise and $\pvariance$ with reducible, parameter uncertainty due to limited data.
We make the common assumption that the measurement noise is mean-zero, independent, identically distributed (IID) and normal $\noise \sim \Nc(\mathbf{0},\sigma^2 \Is)$, where $\sigma$ is, at this point, an unspecified hyperparameter.
In the simplest case, we assume that \emph{both} $\noise$ and $\pvariance$ are normally distributed.
In this case, we need to determine  $\sigma$ and the mean $\bar{\parameters}$ and covariance $\Sigmab_\parameters$ so that
\begin{equation}
\parameters \equiv \bar{\parameters} + \pvariance
\sim \Nc(\bar{\parameters},\Sigmab_\parameters)
\end{equation}
by the well-known sum of Gaussians rule.
The mean parameters $\bar{\parameters}$ and covariance $\Sigmab_\parameters$ are subject to \emph{epistemic convergence}; while the push-forward model response maintains a finite variance.%
\footnote{Epistemic convergence is the process by which model parameters become more certain with increasing amounts of calibration data. By definition, epistemic uncertainty is reducible with more data.}

This formulation naturally leads to Bayesian calibration \cite{kennedy2001bayesian,ghanem2017handbook}.
In fact, the least squares optimization \eqref{eq:least_squares} can be related to Bayesian calibration through Bayes rule
\begin{equation}
\prob(\parameters | \dataset ) = \frac{\prob(\dataset | \parameters)}{\prob(\dataset)} \prob(\parameters)
\end{equation}
with
posterior $\prob(\parameters | \dataset )$,
likelihood $\prob(\dataset | \parameters)$,
prior $\prob(\parameters)$,
and normalizing factor (evidence) $p(\dataset)$.
In the case of maximum likelihood estimation (MLE) (or maximum {\it a posteriori} estimation (MAP) with an uninformative prior $\prob(\parameters) \sim \text{constant}$)
\begin{equation}
\prob(\parameters | \dataset ) = \frac{p(\dataset | \parameters)}{p(\dataset)} \propto \prob(\dataset | \parameters) \equiv \likelihood(\parameters) \ .
\end{equation}
A likelihood $\likelihood$ that corresponds to the MSE objective, $\objective$:
\begin{equation} \label{eq:gaussian_likelihood}
\likelihood(\parameters)
=
\Nc(\outputs; \model(\inputs;\parameters), \Sigmab_\outputs)
= \frac{1}{\sqrt{(2\pi)^{n_\outputs n_r} (\det\Sigmab_\outputs)^{n_r} }}
\exp\left(
-\frac{1}{2} \tr \bigl(\residuals(\inputs;\parameters)^T \Sigmab_\outputs^{-1} \residuals(\inputs;\parameters)\bigr)
\right)
\end{equation}
provides the connection.
Note, $\tr(\residuals \Sigmab^{-1} \residuals^T)$ is a scalar, and units of likelihood, being similar to a probability density, are the inverse of those of the output $\outputs$.
Here again, the expected precision $\Sigmab_\outputs^{-1}$ appears and  the likelihood is a multivariate normal in terms of the residuals $\residuals$, defined in  \eref{eq:residuals}, with parameters $\parameters$ and $\Sigmab_\parameters$.
Hence, the log likelihood is
\begin{eqnarray} \label{eq:log_likelihood}
-2 \log \likelihood(\parameters) &=&
n_\outputs n_r \log( 2\pi) + n_r \log( \det \Sigmab_\outputs)
+ n_r \frac{1}{2} \langle \tr \residuals^T \Sigmab_\outputs^{-1} \residuals \rangle
\end{eqnarray}
where the first two terms do not depend on the parameters $\parameters$ and the last term is the weighted SSE $\Phi$.
Note a double sum is implied in the last term, inner sum over $n_y$ and outer sum over $n_r$, likewise for \eref{eq:gaussian_likelihood}.
and we have introduced the notation
\begin{equation} \label{eq:SSE}
\Phi(\parameters; \dataset, \model)  \equiv
\frac{1}{2} \tr \residuals^T \Sigmab_\outputs^{-1} \residuals
=
\frac{1}{2} n_r \langle \tr \residuals^T \Sigmab_\outputs^{-1} \residuals \rangle
\end{equation}
for an average over the repeats.
This likelihood is consistent with the Kennedy and O'Hagan additive random discrepancy term \cite{kennedy2001bayesian}.

\paragraph{Laplace approximation}
In general, obtaining the moments that determine a Gaussian distribution in \eref{eq:gaussian_likelihood} analytically is difficult because the model $\model$ is non-linear in its parameters $\parameters$.
Laplace approximation is based on the connection between the SSE $\objective$ and the likelihood $\likelihood$ through a Taylor expansion of the log likelihood \eqref{eq:log_likelihood}:
\begin{equation}
\objective \approx \objective(\parameters^*)
+ \underbrace{\partialb_\parameters \objective(\parameters^*) }_{\mathbf{0}} \Delta\parameters
+ \frac{1}{2} \underbrace{ \partialb_\parameters \partialb_\parameters \objective(\parameters^*)}_{\hessian} \Delta\parameters \otimes \Delta\parameters
\end{equation}
where $\Delta\parameters \equiv \parameters - \parameters^*$ is a perturbation around $\parameters^*$ and the second term is zero at an optimum of $\Phi$, as defined in \eref{eq:least_squares}.
The Hessian $\hessian$ at $\parameters = \parameters^*$ consists of two terms:
\begin{equation}
\hessian \equiv \partialb_\parameters \partialb_\parameters \objective
= \frac{1}{2} \partialb_\parameters  ( \partialb_\parameters \model \cdot \Ws \residuals )
=  \frac{1}{2} \partialb_\parameters \model \otimes \Ws \partialb_\parameters \model
+ \frac{1}{2} \partialb^2_\parameters \model \cdot \Ws \residuals
\approx \frac{n_r}{2}  \sensitivities^T \Ws \sensitivities \ ,
\end{equation}
the first is formed from the parameter sensitivities of the model $\sensitivities\equiv\partialb_\parameters \model$, and the second term can be omitted if the model discrepancies are negligible $\residuals \approx \zerob$ or the model is linear in its parameters.%
\footnote{Note the Gauss-Newton approximation $\hessian \approx \tfrac{n_r}{2}  \sensitivities^T \Ws \sensitivities $ can be corrected by a factorization of the symmetric second order term $\Bs^T\Ws\Bs =\partialb^2_\parameters \model \cdot \Ws \residuals$ so that $\As \to \As + \Bs$.
Also, the Gauss-Newton approximation of the Hessian is essentially equivalent to linearizing the model then forming the Hessian, as in linear Bayesian optimization.}
Here $\sensitivities$ are the (mean) parameter sensitivities at $\parameters=\parameters^*$:
\begin{equation}\label{eq:sensitivities}
\sensitivities(\inputs) \equiv
\langle \partialb_\parameters \residuals \rangle
= \partialb_\parameters \model(\inputs_0; \parameters^*) \ .
\end{equation}
where $\inputs_0$ is the input for the first repeat (which is equal to all other repeats).
Now using the association of $\log \prob$ with $\Phi$ from \eref{eq:SSE} we can form the approximation
\begin{equation}
\likelihood \equiv
\prob(\parameters | \dataset) \approx
\sqrt{\frac{\det \hessian}{(2\pi)^{n_\parameters}}}
\exp\left( - \frac{1}{2} (\parameters - \parameters^*) \Hs (\parameters-\parameters^*) \right)
\end{equation}
where the Hessian acts as a precision and we have arbitrarily set $\objective(\parameters^*) = 0$.

\paragraph{Embedded Laplace Method}
Here we enhance the classical Laplace approximation by adopting an embedded uncertainty \cite{sargsyan2015statistical} perspective where we treat the model parameters as random variables, e.g.
\begin{equation}
\parameters \equiv \bar{\parameters} + \pvariance = \bar{\parameters} + \varthetab \zetab
\end{equation}
with mean $\bar{\parameters} = \EXP[\parameters]$, standard deviation $\varthetab = \sqrt{\VAR[\parameters]} \equiv \sqrt{\Sigmab_\parameters}$, expressed as a Cholesky decomposition $\Ls \Ls^T = \varSigmab_\parameters$, and $\zetab \sim \Nc(\mathsf{0},\mathsf{I})$, and aim to determine the moments of that define these distributions.
In this context, the parameters defining the distributions converge with increasing data, while distributions themselves maintain a finite width.

We can obtain the variance of the parameters $\parameters = \bar{\parameters} + \pvariance$ through the linearization of the model
at an optimum $\parameters^*$ in \eref{eq:least_squares}, so that
\begin{equation}
\langle \residuals \rangle \equiv \outputs - \model(\inputs; \parameters^*) \approx
\sensitivities(\inputs) \pvariance + \noise \ .
\end{equation}
If we assume the expectation of $\residuals$ is approximately zero, then
\begin{equation} \label{eq:pcorr}
\Sigmab_\outputs \approx  \sensitivities \Sigmab_\parameters \sensitivities^T  + \sigma^2 \Ib
\end{equation}
via the variance rule of a Gaussian of a linear transformation, $\operatorname{var} [ \sensitivities \inputs ] = \sensitivities \Sigmab_\inputs \sensitivities^T$.%
\footnote{
With more detail:
\begin{equation*}
\Sigmab_\outputs
\equiv \bigl< (\outputs - \langle \outputs \rangle)^2 \bigr>
= \bigl< ( \model(\inputs; \parameters) + \noise -  \langle \outputs \rangle)^2 \big>
\approx \bigl< ( \model(\inputs; \bar{\parameters})  + \sensitivities \pvariance + \noise -  \langle \outputs \rangle)^2 \big>
\approx \bigl< ( \sensitivities \pvariance + \noise )^2 \big>
= \sensitivities \Sigmab_\parameters \sensitivities^T + \sigma^2 \Ib
\end{equation*}
given $\pvariance$ is uncorrelated with $\noise$, $\pvariance \sim \Nc(\mathsf{0},\Sigmab_\parameters)$ and $\noise \sim \Nc(\mathsf{0}, \sigma^2 \mathsf{I})$.}
Note $\sensitivities \Sigmab_\parameters \sensitivities^T$ has the rank of $\Sigmab_\parameters$, i.e. $n_\parameters$, and the dimensions of $\Sigmab_\outputs$, i.e. $n_y \times n_y$, such that $\sensitivities$ is $n_y \times n_\parameters$ and $n_y \gg n_\parameters$.

If we examine \eref{eq:pcorr}, it is clear that $\sensitivities \Sigmab_\parameters \sensitivities^T$ has the rank of $\Sigmab_{\parameters}$ and, hence, at most rank the number of parameters $n_{\parameters}$.
For this work, we will assume $\Sigmab_\parameters$ is full rank due to the experimental loading fully exercising the model.
So the solution for $\Sigmab_\parameters$ and $\sigma$ splits into the span of $\sensitivities \Sigmab_\parameters \sensitivities^T$ and its complement.
Since the complement only depends on $\sigma$, we use it to determine the aleatoric noise scale $\sigma$.
By applying an eigenvalue decomposition of $\Sigmab_\outputs$,
\begin{equation}
\Sigmab_\outputs = \Us \Ds \Us^T
\end{equation}
the largest values determine $\Sigmab_\parameters$ and we can estimate $\sigma^2$ with the average of the lowest singular values.
Since the data can have more information content than the model applied to the same loading we omit the first $\max(n_r,n_\parameters)$ singular values in this estimate.
With $\sigma$ determined, we have a least squares solution for $\Sigmab_\parameters$
\begin{eqnarray} \label{eq:LS_solve}
\sensitivities \Sigmab_\parameters \sensitivities^T &=&  \Sigmab_\outputs-\sigma^2 \Ib \\
\sensitivities^T \sensitivities \Sigmab_\parameters \sensitivities^T \sensitivities  &=& \sensitivities^T \Sigmab_\outputs \sensitivities -  \sigma^2 \sensitivities^T \sensitivities \nonumber \\
\Sigmab_\parameters &=& ( \sensitivities^T \sensitivities )^{-1} \sensitivities^T \Sigmab_\outputs \sensitivities ( \sensitivities^T \sensitivities )^{-1}  -  \sigma^2 (\sensitivities^T \sensitivities)^{-1} \nonumber
\end{eqnarray}
In practice since we already have the eigensystem, we use the largest $n_\parameters$ eigenvectors to form  $\check{\Us} = \Us_{1:n_\outputs,1:n_\parameters}$ and $\check{\Ds} = \Ds_{1:n_\parameters,1:n_\parameters} - \sigma^2$ and then solve
\begin{equation} \label{eq:solve_Sigma}
\Sigmab_\parameters = (\check{\Us}^T \sensitivities)^{-1} \check{\Ds} (\check{\Us}^T \sensitivities)^{-T} \ .
\end{equation}
where, for numerical conditioning we use the (Moore-Penrose) pseudoinverse for $(\Us^T \sensitivities)^{-1}$.

\paragraph{Fit}
Subsequently, we can fit the parameter covariance using an optimization objective.
In fact, we can determine $\Sigmab$ by optimizing log posterior predictive \eqref{eq:log_likelihood} with the linearization of the model response
\begin{equation} \label{eq:max_lop_posterior}
\Sigmab_\parameters, \sigma=
\Argmin{\Sigmab_\parameters \succeq 0, \sigma \ge 0 } \left( -\log \prob(\parameters | \dataset) \right)  =
\Argmin{\Sigmab_\parameters \succeq 0, \sigma \ge 0 } \left( n_r \log( \det \Sigmab_\outputs)
+ \tr ( \residuals  \Sigmab_\outputs^{-1} \residuals) \right)
\end{equation}
given $\parameters = \parameters^*$.
This optimization minimizes the SSE (the second term) regularized by the first term.

To start the optimization we use \eref{eq:solve_Sigma} with a sample covariance for $\Sigmab_\outputs \approx \tfrac{1}{n_r-1} \Rs \Rs^T$  as an initial guess.
Also, for numerical conditioning, we decompose the covariance matrix $\Sigmab_\outputs$ into variances and correlation coefficients \cite{fisher1915frequency,fisher1921probable,archakov2021new}:
\begin{equation}
\left[ \Sigmab_\outputs \right]_{ij} = \sqrt{\varsigma^2_{i}}\left[\boldsymbol{\Xi}\right]_{ij} \sqrt{\varsigma^2_{j}}
\end{equation}
with the diagonal of the covariance matrix $\varsigma_i^2 = [ \Sigmab_\outputs ]_{ii}$ and the correlation $\Xi_{ij}$ with off-diagonal components $\rho_{ij}$, which are reparameterized as: $\varsigma^2 \leftarrow \exp(\sigma^2)$ and $\rho \leftarrow \arctanh(\rho)$.

\paragraph{Linear-in-parameters models}

The method is exact for linear-in-parameters models, and this case provides an analytical illustration.
A general form for a model that is linear in its parameters $\parameters$ is:
\begin{equation}
\model(\inputs) = \sum_a \phi_a(\inputs) \theta_a = \sensitivities(\inputs) \parameters \ ,
\end{equation}
where the sensitivities are composed of basis-like functions $\sensitivities(\inputs) = [ \phi_a(\inputs) ]_a$;
for instance, for a model that is also linear in its inputs $\inputs$
\begin{equation}
\sensitivities
= \begin{bmatrix} \mathsf{1}, \mathsf{1}, ... \mathsf{1} \\ \inpt_0, \inpt_1, ... \inpt_n \end{bmatrix}^T
=  [ \mathsf{1}, \inputs ]^T \ .
\end{equation}
where, here, we use $\inputs$ for its average over repeats $\langle \inputs \rangle$ for clarity.
In this case, \eref{eq:embedded_model} becomes:
\begin{equation}
\outputs \approx \sensitivities(\inputs)\Thetab + \Ns  \quad \Rightarrow \quad
\residuals \equiv \outputs - \sensitivities(\inputs)\Thetab - \Ns
\end{equation}
where $\Ns$ are samples of $\noise \sim \Nc(0,\sigma^2)$ and $\Thetab$ are samples of $\thetab \sim \Nc(\bar{\thetab}, \Sigmab_\parameters)$,
And the optimum $\thetab^* = \bar{\thetab} \equiv (\sensitivities^T \sensitivities)^{-1} \sensitivities^T \langle{\outputs}\rangle$.
Hence,
\begin{equation}
\EXP[\residuals] = \EXP[\outputs] - \sensitivities(\inputs)\EXP[\Thetab] - \EXP[\Ns]
= \EXP[\outputs] - \sensitivities(\inputs) \bar{\thetab}
\end{equation}
and
\begin{eqnarray}
\Sigmab_\outputs &=& \VAR[\outputs\outputs^T] = \EXP[\outputs\outputs^T] - \EXP[\outputs] \, \EXP[\outputs] \\
&=& \sensitivities \VAR[\Thetab \Thetab^T ] \sensitivities^T  + \EXP[\Ns \Ns^T] + \sensitivities \EXP[\Thetab \Ns^T ] + \EXP[\Ns \Thetab^T ] \sensitivities^T \nonumber \\
&=&
\sensitivities \Sigmab_\parameters \sensitivities^T + \sigma^2 \Is \nonumber
\end{eqnarray}
if we  assume: (a)
$\EXP[\Ns \Thetab^T ] = \EXP[\Ns] \EXP[\Thetab^T] + \VAR[\Ns, \Thetab^T] = \zerob$
and (b)
$\EXP[\residuals] = \zerob$ so that $\EXP[\outputs] = \sensitivities(\inputs) \bar{\thetab}$.

\paragraph{Variational inference}
Like the embedded Laplace fitting technique,
the VI technique minimizes the difference between the posterior and a surrogate $q_\phib (\parameters)$ through optimization of the parameters $\phib$.
The primary objective is the Kullback-Leibler divergence (KLdiv):
\begin{equation}
\KLD(q_{\phib}(\parameters) || p(\parameters | \dataset)) =
\EXP_{q_{\phib}}\left[ \log \frac{q_{\phib}(\parameters)}{p(\parameters | \dataset)} \right] \ ,
\end{equation}
but the evidence lower bound (ELBO) is a more tractable optimization objective:
\begin{equation}
\operatorname{ELBO}(q_{\phib}(\parameters) || p(\parameters | \dataset))
= \EXP_{q_\phi}\left[ \log p(\dataset | \parameters) \right] - \KLD( {q_\phi} (\parameters) \| p(\parameters) )
\end{equation}
that combines a data consistency objective with a prior similarity regularizing term.

In the \emph{black box variational inference} (BBVI) \cite{ranganath2014black} formulation of \emph{mean field} VI a Gaussian with a diagonal covariance (with elements $\sigma_i$) and a uniform prior allow the reduction of the ELBO to:
\begin{equation} \label{eq:indep_elbo}
\operatorname{ELBO}(q_{\phib}(\parameters) || p(\parameters | \dataset))
= \log \likelihood(\parameters) +    \underbrace{\frac{n_\phi}{2}(1 + \log 2\pi ) + \sum_i \log \phi_i}_{\entropy(q_\phi)}
\end{equation}
The first term in \eref{eq:indep_elbo} is the log likelihood from \eqref{eq:log_likelihood} (ignoring constant terms) and is model dependent, while the remainder only depends on the Gaussian prior, where $\entropy$ is the entropy.
For comparison, the entropy of full covariance Gaussian is:
\begin{equation}
H(q_\phi) = -\EXP\left(\log q_\phi\right) = \frac{n_\phi}{2} \left( 1 + \log(2\pi) \right) + \frac{1}{2} \log \det \Sigmab
\end{equation}

\paragraph{Embedded variational inference}

To accommodate the embedded error formulation \eqref{eq:embedded_model} within VI, we take the mean parameters from the MAP and determine the variance parameters through (subsequent) optimization.
In particular we use the form:
\begin{equation}
\parameters = \bar{\parameters} + \varparameters \odot \zetab
\end{equation}
where $\bar{\parameters}$ are the mean parameters and $\varparameters$ scale the contributions of independent Gaussian random variables $\zeta_i \sim \Nc(0,1)$.
This is consistent with a mean field approximation and, hence, is in contrast with a  Cholesky-like structure that would allow cross-correlations.
For this formulation we chose a likelihood consistent with the \emph{Gaussian approximation to marginalized likelihood}
from \cref{sargsyan2015statistical}:
\begin{equation}
\likelihood(\parameters; \data) \approx \prod_i \Nc(\outpt_i;  \mub(\inpt_i; \parameters); \Sigmab(\inpt_i; \parameters))
\end{equation}

\paragraph{Error metrics}
To assess deterministic quantities we use the (usual) root mean square error (RMSE):
\begin{equation}
\operatorname{RMSE}(\parameters,\parameters^*) = \sqrt{ \frac{1}{N} \sum_{i=1}^{N} (\theta_i - \theta_i^*)^2}
\end{equation}

To assess probabilistic predictions, we use the Wasserstein-1 ($\Wc_1$) distance to compare the associated distributions:
\begin{equation}
\mathcal{W}_1(\prob_a, \prob_b)
= \int \vert \CDF_a(\outpt') - \CDF_b(\outpt')\vert \, \mathrm{d}\outpt' \ ,
\end{equation}
where the cumulative distribution functions (CDFs) result from empirical probability density functions (PDFs):
\begin{equation}
\CDF(\outpt) = \int_{-\infty}^{\outpt} \prob(\outpt') \, \mathrm{d}\outpt'
= \int_{-\infty}^{\outpt}  \sum_i \delta(\outpt'-\outpt_i) \, \mathrm{d}\outpt'
\end{equation}
For two univariate Gaussians with the same mean, the Wasserstein-1 reduces to:
\begin{equation} \label{eq:sigma_scale}
\Wc_1 =  \frac{2}{\sqrt{2\pi}} \vert \sigma_2 - \sigma_1 \vert
\end{equation}
which provides the expected scaling of $\mathcal{W}_1$ with the standard deviations of the PDFs being compared.
For multivariate Gaussians, $\Nc(\mub_1,\Sigmab_1)$ and  $\Nc(\mub_2,\Sigmab_2)$, $\Wc_1$ can be bounded:
\begin{equation}
\| \mub_2-\mub_1 \|_2  \le \mathcal{W}_1 \le \sqrt{\Wc_2} = \sqrt{
\| \mub_1 - \mub_2\|_2^2 + \tr (\Sigmab_1 + \Sigmab_2) - 2 \tr (\sqrt{\Sigmab_1\Sigmab_2} )
} \ .
\end{equation}
which provides additional intuition.

\section{Demonstrations} \label{sec:results}
First we use simple calibrations to illustrate epistemic and aleatory behavior of the embedded Laplace (EL) and compare it with BBVI.
Then, we apply the EL method to real-world calibration of complex material models to experimental data.

\subsection{Illustrations}

For the illustrations, we employ the two parameter models listed in \tref{tab:models} for easy interpretation and visualization.
The linear and quadratic models are linear in their parameters and the third is nonlinear.
For each data generating model, we take realizations of the parameters $(a,b)$ from $\Nc(\mub,\Sigmab)$ where:
\begin{equation}
\mub =
\begin{bmatrix}
2.0  \\
1.0 \\
\end{bmatrix}
\quad
\Sigmab =
\alpha
\begin{bmatrix}
0.4  & 0.05 \\
0.05 & 0.1 \\
\end{bmatrix}
\end{equation}
and $\alpha$ is a scaling factor.
The standard deviation of the external Gaussian noise $\noise$ was set to $\sigma = 0.01$.

\begin{table}[h!]
\centering
\begin{tabular}{|l|l|}
\hline
linear    & $a x   + b$  \\
quadratic & $a x^2 + b x$  \\
\hline
nonlinear &  $a (1-\exp(-b x))$ \\
\hline
\end{tabular}
\caption{Two parameter models}
\label{tab:models}
\end{table}

\subsubsection{Aleatory/epistemic behavior}
\begin{figure}
\centering
\begin{subfigure}[c]{0.32\linewidth}
\includegraphics[width=0.99\linewidth]{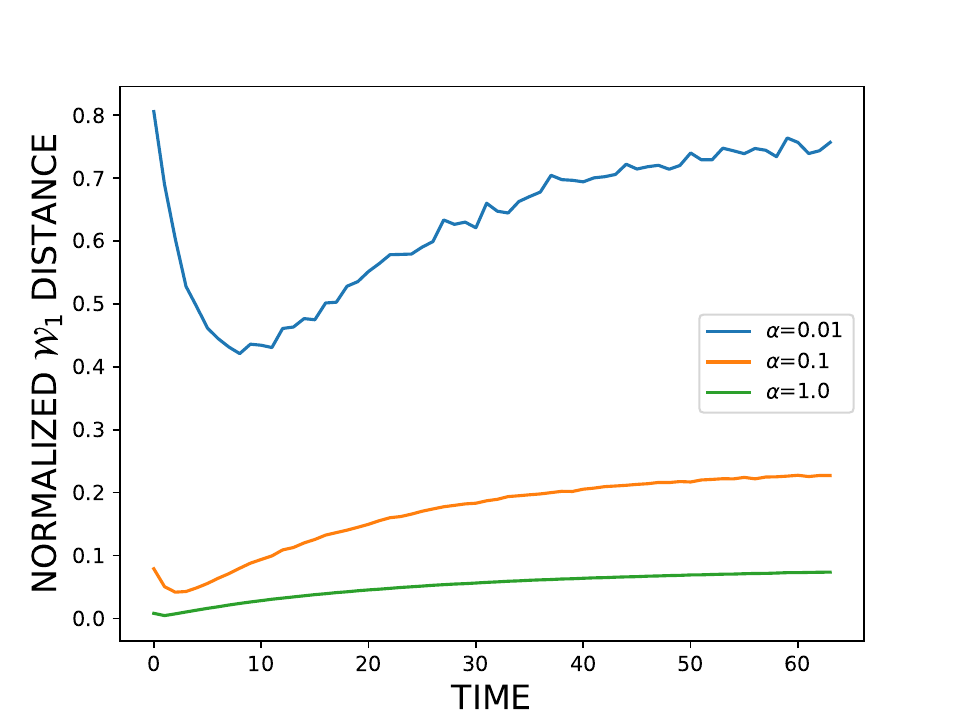}
\caption{$\Wc_1$ distance vs time}
\end{subfigure}
\begin{subfigure}[c]{0.32\linewidth}
\includegraphics[width=0.99\linewidth]{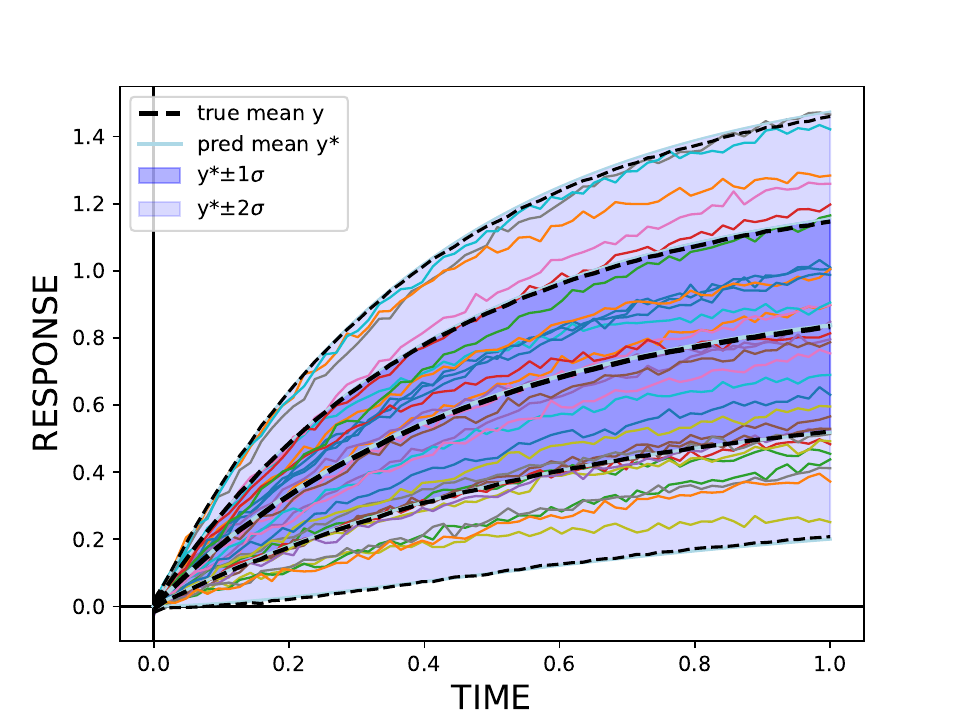}
\caption{response $\outpt(t)$ for large variance}
\end{subfigure}
\begin{subfigure}[c]{0.32\linewidth}
\includegraphics[width=0.99\linewidth]{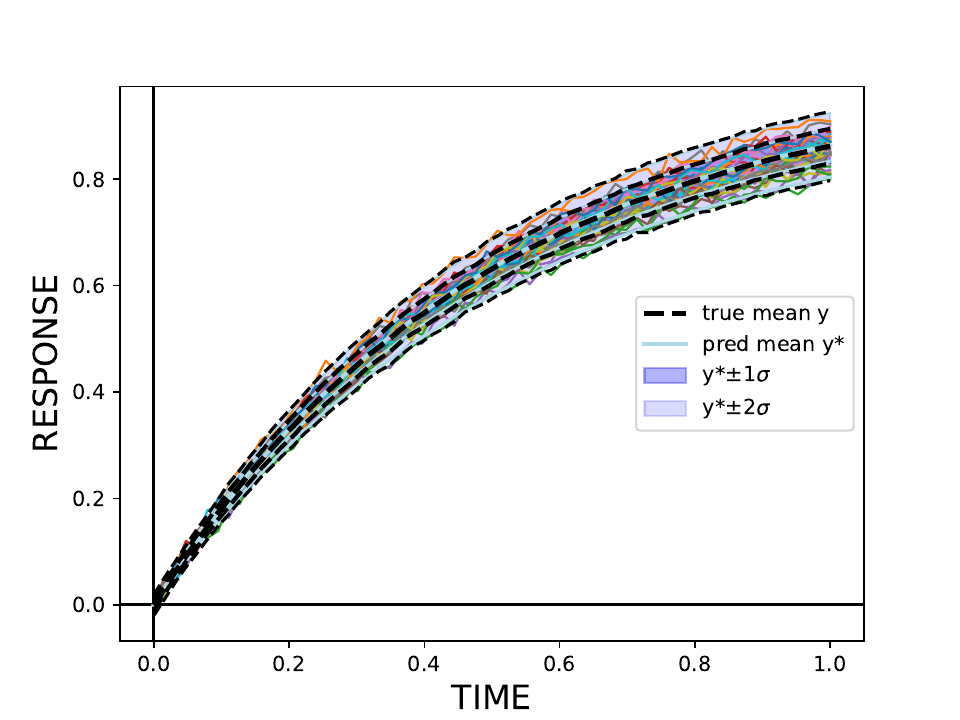}
\caption{response $\outpt(t)$ for small variance}
\end{subfigure}
\caption{Illustration: (a) convergence in terms of a normalized $\Wc_1$ distance and prediction comparisons for  (b) high parametric variance, $\alpha=1.0$ (c) low  parametric variance, $\alpha=0.01$. }
\label{fig:illustration_variance}
\end{figure}

In this study, we fix the measurement noise and change parametric variance via $\alpha$.
\fref{fig:illustration_variance}a shows the relative error as measured by $\Wc_1/\sigma$ (as motivated by \eref{eq:sigma_scale}.
The trend is for the relative error to improve as the parametric variance becomes larger than the noise.
Clearly, \fref{fig:illustration_variance}b,c show that the EL method provides predictions, as indicated by the 2 standard deviation bands, that cover the data, as depicted by the trajectory realizations, in both the large and small variance cases.

\subsubsection{Comparison with VI}

In this illustration, we compare BBVI to the proposed EL method.
For the BBVI approach, we assume independent internal variance
\begin{equation}
\parameters = \bar{\parameters} + \varparameters \odot \zetab, \quad \zetab \sim \Nc(\mathbf{0},\mathbf{I})
\end{equation}
in analog with mean field VI.
Specifically
\begin{equation}
\theta_i = \bar{\theta}_i + \vartheta_i \zeta, \quad  \zeta \sim \Nc(0,1), \quad i \in (1,2)
\end{equation}
so that 4 parameters ($\bar{\theta}_1, \vartheta_1, \bar{\theta}_2, \vartheta_2$) need to be determined.

\begin{figure}
\centering
\begin{subfigure}[c]{0.32\linewidth}
\includegraphics[width=0.99\linewidth]{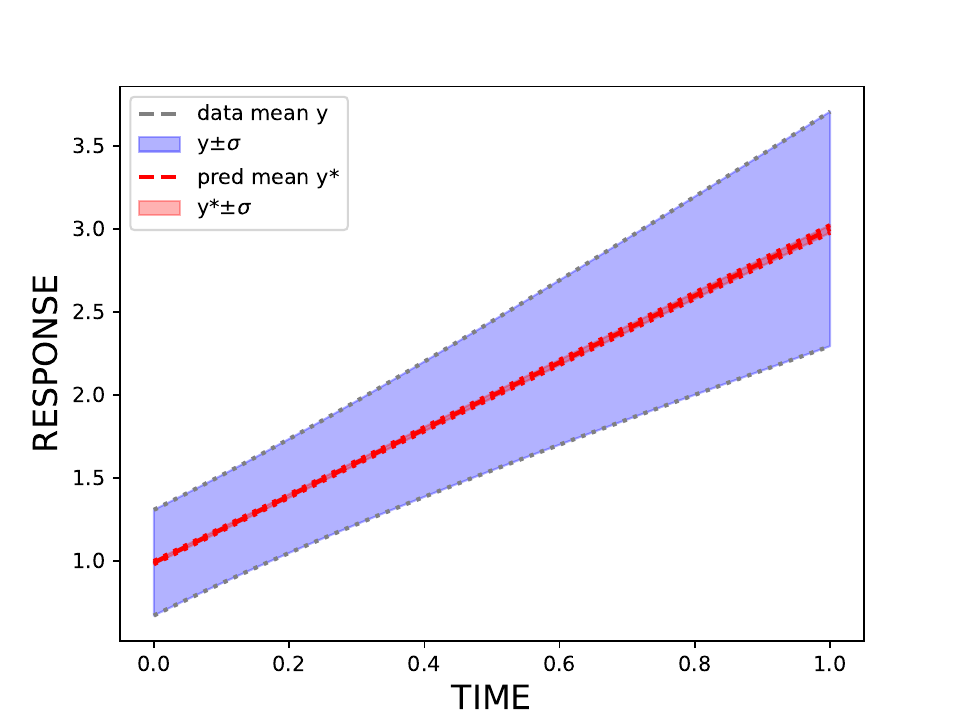}
\caption{VI}
\end{subfigure}
\begin{subfigure}[c]{0.32\linewidth}
\includegraphics[width=0.99\linewidth]{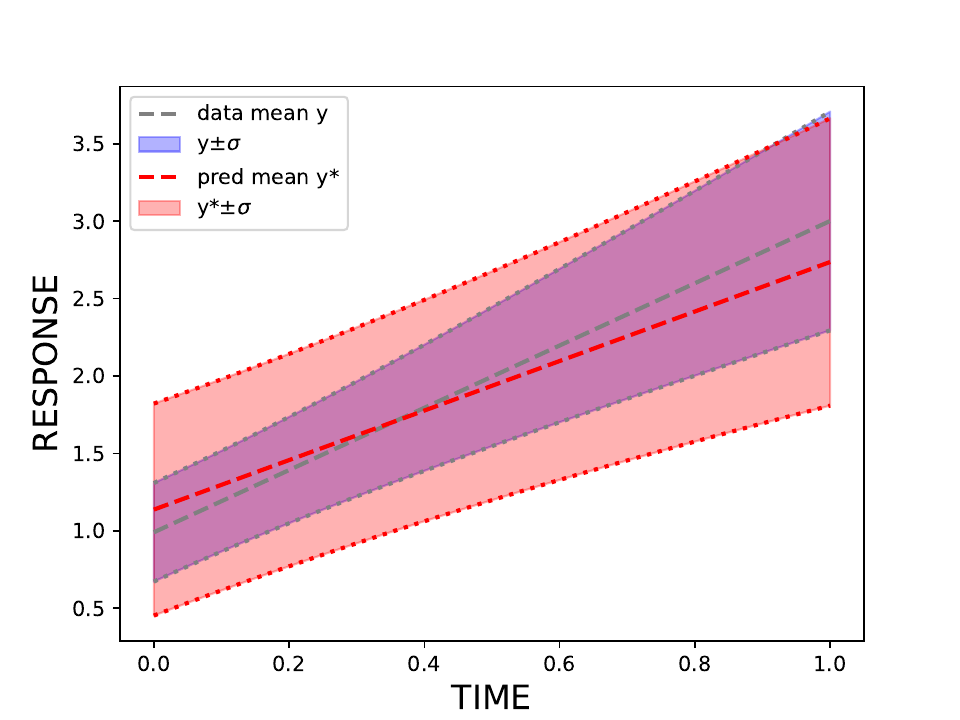}
\caption{embedded VI}
\end{subfigure}
\begin{subfigure}[c]{0.32\linewidth}
\includegraphics[width=0.99\linewidth]{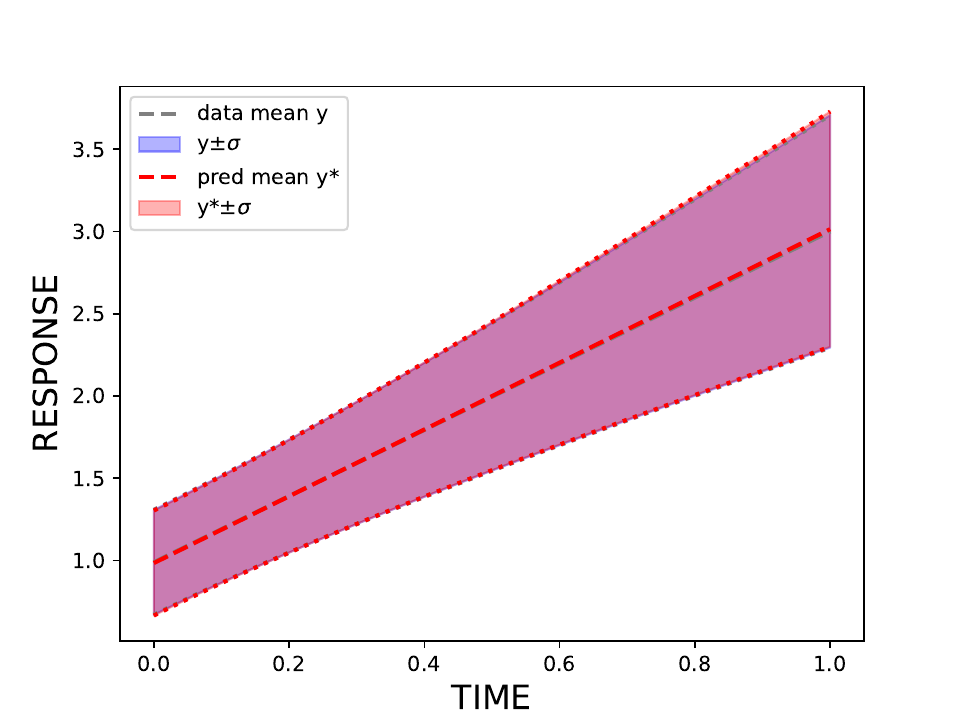}
\caption{EL}
\end{subfigure}
\caption{Illustration: comparison of VI, embedded VI, EL push-forwards for a linear model with parametric variance (blue: data, red: prediction).}
\label{fig:illustration_comparison}
\end{figure}

\begin{figure}
\centering
\includegraphics[width=0.45\linewidth]{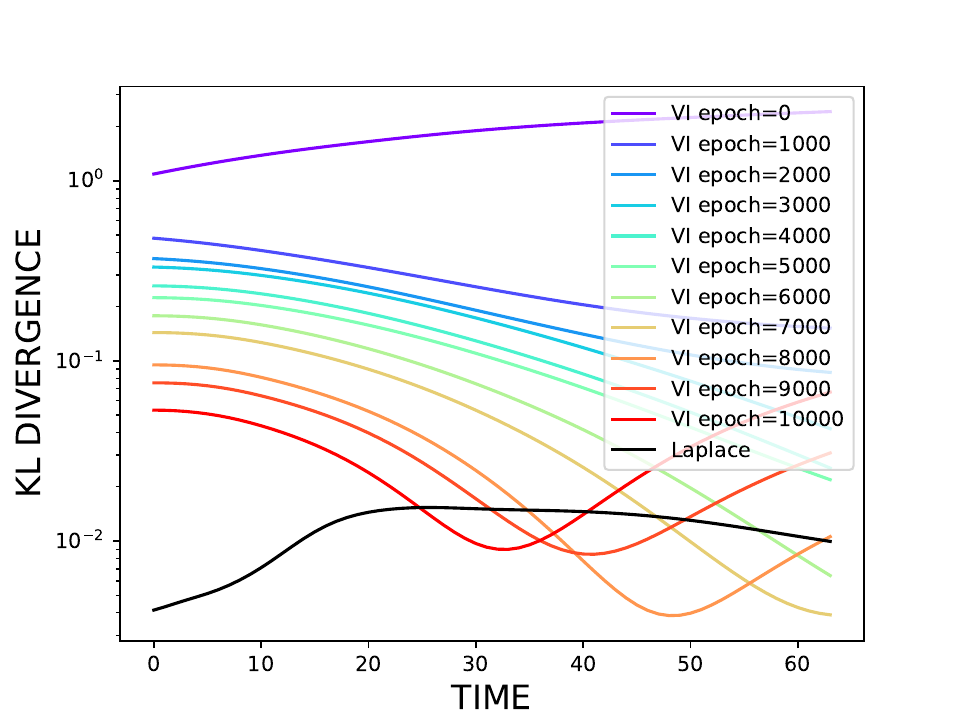}
\caption{Illustration: convergence of embedded VI compared to EL in terms of the KL divergence between the underlying data variance and the predicted posterior.
Note the difference in the KL divergence profiles for EL and the final embedded VI in \fref{fig:illustration_convergence} are due to the assumed covariance structure, i.e. full for EL and diagonal for the embedded VI result.
}
\label{fig:illustration_convergence}
\end{figure}

For the EL method, we assume we can afford the cross-terms typically omitted in (\emph{mean field}) VI.
\fref{fig:illustration_comparison} plots convergence of both the standard mean field VI and the embedded version versus the EL result.
With standard VI the results collapses to the mean, as expected, while mean-field embedded version of VI maintains a finite variance but has errors due to the diagonal covariance assumption.
On the other hand, the EL with a full covariance captures the data distribution accurately.
With reference to \fref{fig:illustration_convergence}, clearly it takes thousands of epochs of gradient descent with a commensurate number of function evaluations.
The optimization method (e.g. switching to a quasi-Newton method like BFGS \cite{nocedal2006numerical}) or learning rate (0.01 for this illustration) will change the convergence but not the conclusion that EL is substantially more model sample efficient.
For this example, EL is both more accurate in terms of estimating the true variance and more efficient in terms of model evaluations than standard VI and a version that accommodates embedded parameter error.

\subsubsection{Model discrepancy}

\begin{figure}
\centering
\begin{subfigure}[c]{0.65\linewidth}
\includegraphics[width=0.95\linewidth]{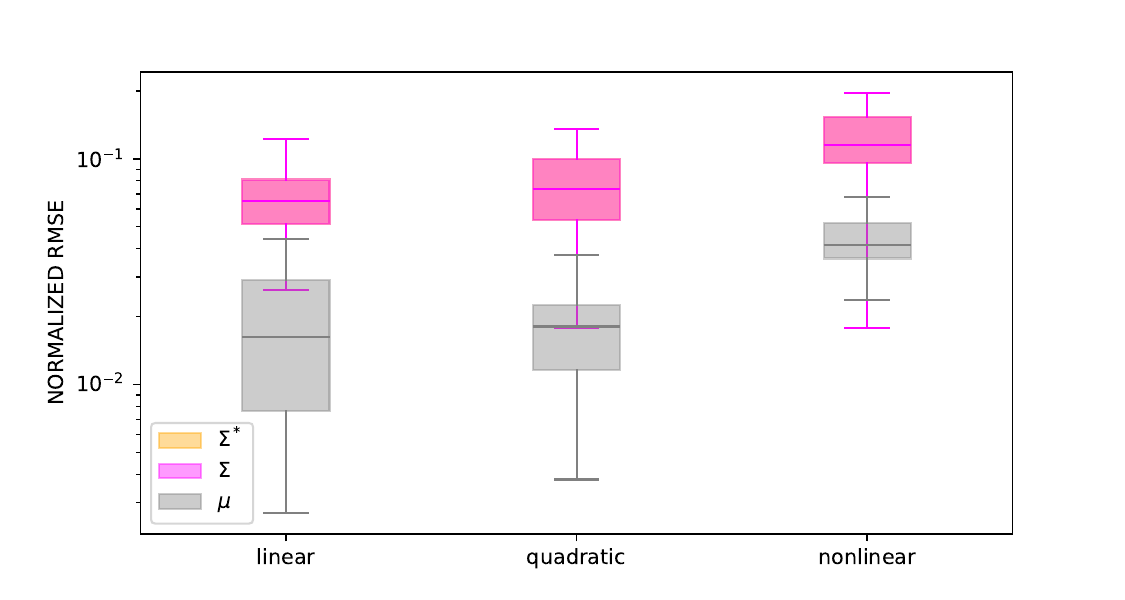}
\caption{parameter}
\end{subfigure}
\begin{subfigure}[c]{0.65\linewidth}
\includegraphics[width=0.95\linewidth]{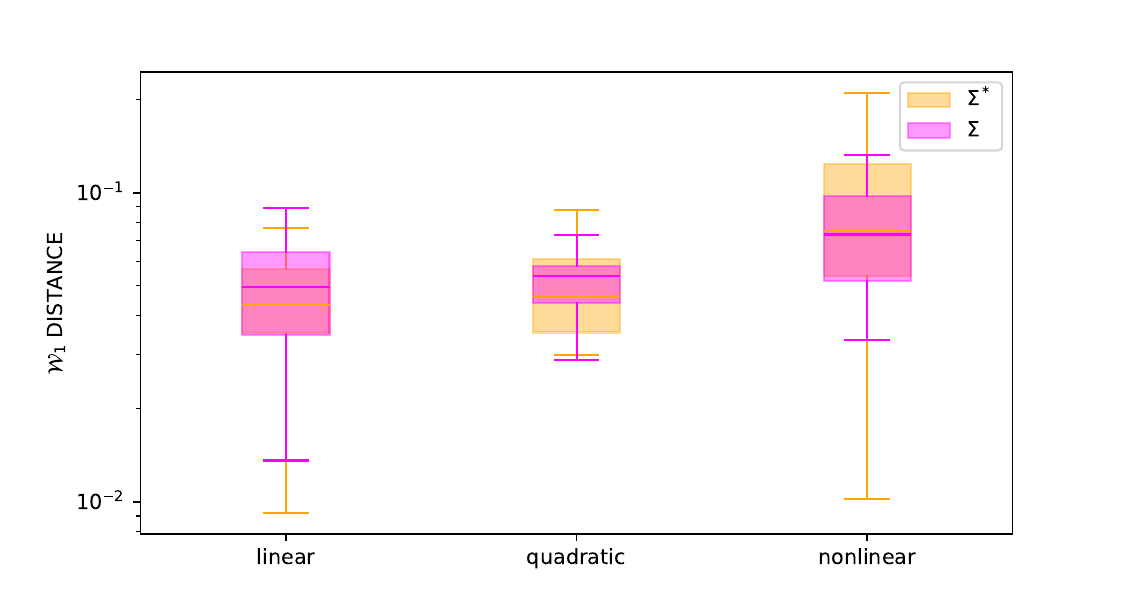}
\caption{response}
\end{subfigure}
\caption{Illustration: RMSE of the estimated mean and variance for no model discrepancy calibrations (a) parameter errors, (b) push-forward errors.}
\label{fig:no_discrepancy}
\end{figure}

\begin{figure}
\centering
\begin{subfigure}[c]{0.65\linewidth}
\includegraphics[width=0.95\linewidth]{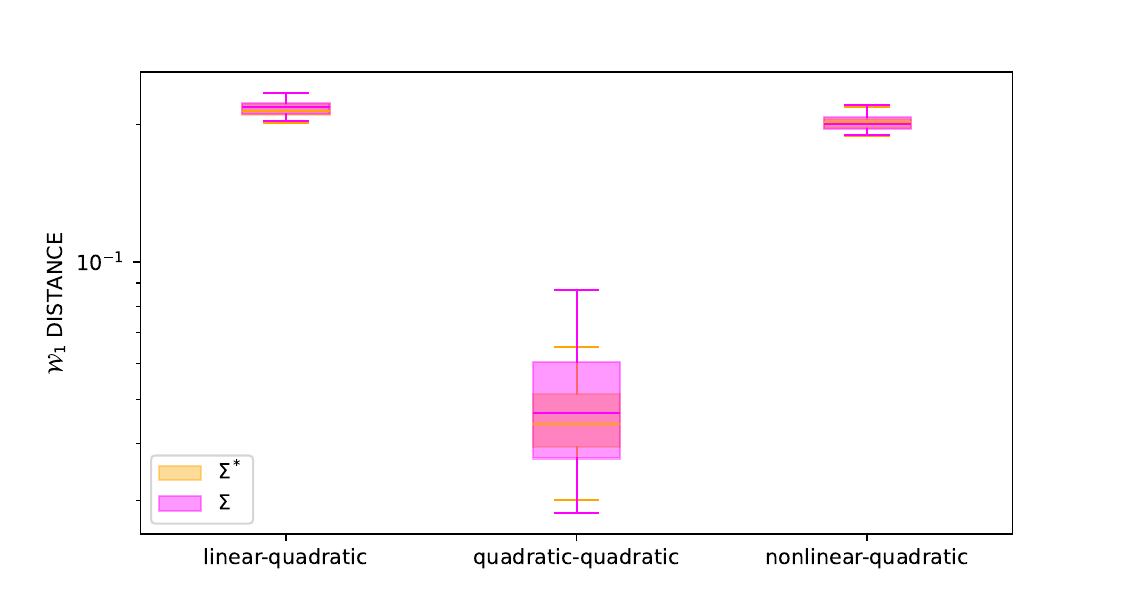}
\caption{quadratic data}
\end{subfigure}
\begin{subfigure}[c]{0.65\linewidth}
\includegraphics[width=0.95\linewidth]{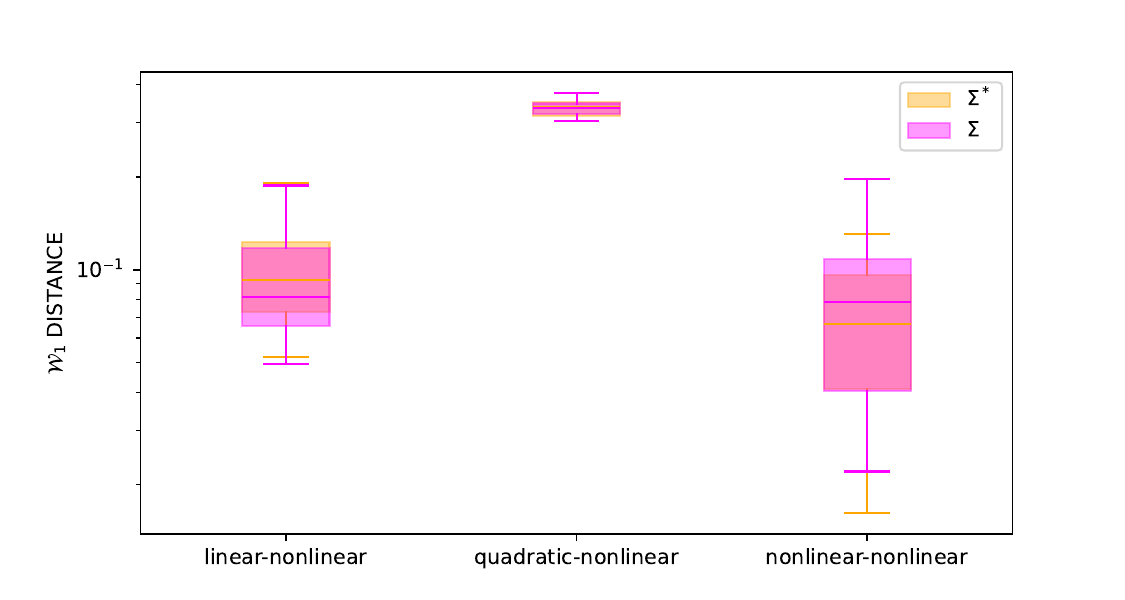}
\caption{nonlinear data}
\end{subfigure}
\caption{Illustration: $\Wc_1$ distance for model discrepancy (a) quadratic data, (b) nonlinear data.}
\label{fig:discrepancy}
\end{figure}

In this last example we illustrate the behavior of the EL method with significant model discrepancy.
For reference,  \fref{fig:no_discrepancy} shows the errors for the case where the data-generating and the fitting model are the same.
The parameter means have lower errors than the covariance (a higher statistical moment), while errors for the linear algebra guess $\Sigmab_\parameters^*$ are on par with the calibrated $\Sigmab_\parameters$.
On the other hand, the push-forward response statistics show that the fitted covariance have marginally better error in some cases.
\fref{fig:discrepancy} shows corresponding results for the fitting model being not necessarily the same as the data generating model.
In this case, if the fitting model is sufficiently similar to the data generating model,e.g. for a quadratic model fitted to nonlinear data, the errors, as measured by $\Wc_1$, are on par with the nonlinear model fitted to the same data; otherwise, the $\Wc_1$ distance is substantially larger.
Nevertheless, these results  also show that the fitted $\Sigmab_\outputs$  can improve on the estimate $\Sigmab_\outputs^*$ but is not necessary for an engineering-level result.

\subsection{Applications}
Our engineering-oriented applications involve experimental data and empirical material models.
In particular, many high-consequence performance evaluations require predicting the likelihood of failure due to plastic deformation.

Plasticity is a history-dependent material response that has a recoverable elastic regime and an irreversible plastic regime that evolves with loading.
The basic framework of these models is: (a) a stress rule
\begin{equation}
\stress =  \stress(\strain,\internalstate)
\end{equation}
relating stress $\stress$ to strain $\strain$ and internal state $\internalstate$, (b) a yield condition
\begin{equation}
\yield = \| \stress \| - h(\internalstate) \le 0
\end{equation}
that delineates elastic ($\yield < 0$) from plastic ($\yield = 0$) behavior, and (c) a flow rule
\begin{equation}
\dot{\internalstate} = \fb(\strain,\internalstate)
\end{equation}
that governs the evolution of internal state $\internalstate$ when $\yield = 0$.
Together the components of (a-c) elastoplasticity constitute a complex differential-algebraic equation (DAE) system \cite{gear1984differential,ascher1998computer}.
More details can be found in texts \cite{lubliner2008plasticity,gurtin2010mechanics} and the papers \cite{karlson2023toward} these calibrations are based on.
All the uncertain parameters are associated with the Voce-like \cite{voce1948relationship,voce1955practical,swift1952plastic,johnson1983constitutive} hardening rule:
\begin{equation} \label{eq:voce}
h(\epsilon_p,\dot{\epsilon}_p) = Y (1+ c\log(\dot{\epsilon}_p/\dot{\epsilon}_0)) + A (1 - \exp(-b \epsilon_p)) \ ,
\end{equation}
that controls the radius of the yield surface $\yield=0$, and is dependent on  the equivalent plastic strain $\epsilon_p$ (an evolving internal state variable) and its rate.
Setting $c=0$ removes the rate dependence and reduces the model to the standard Voce form.
Also, the relative sensitivity to $b$ is typically high due to it appearing inside the exponential function.
We assume the elastic parameters are precisely determined by the baseline calibration and do not vary across the test specimens.

For these demonstrations the model sensitivities around the MLE fit $\parameters^*$ are generated with finite differences.
The parameter sensitivities are gradients which require $\propto n$ function evaluations, where the exact coefficient depends on the scheme, likewise a full Hessian $\propto n^2 $ evaluation.
We use the noise estimate $\sigma < \Delta X$  to set the finite difference step size $\Delta X \approx | \Delta \inputs |$.

\subsubsection{Synthetic validation}
\begin{figure}
\centering
\includegraphics[width=0.30\linewidth]{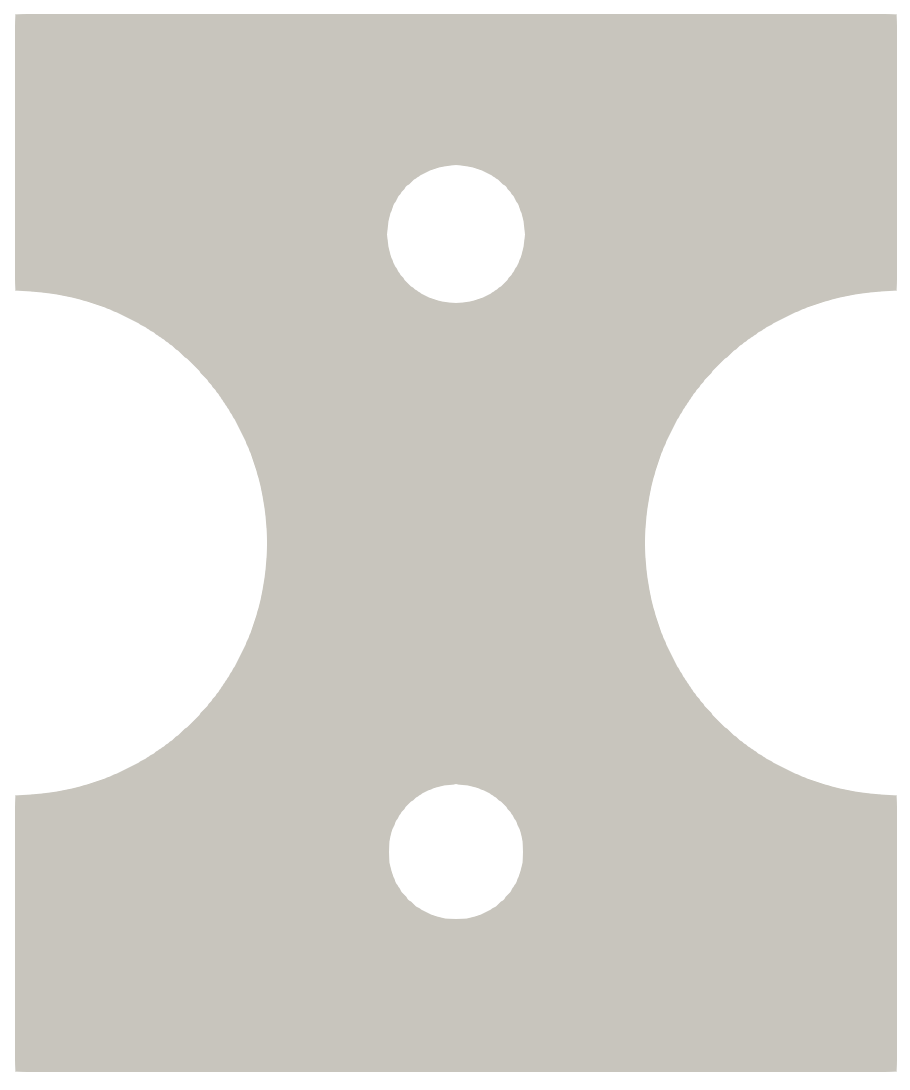}%
\hspace{0.1\linewidth}%
\includegraphics[width=0.25\linewidth]{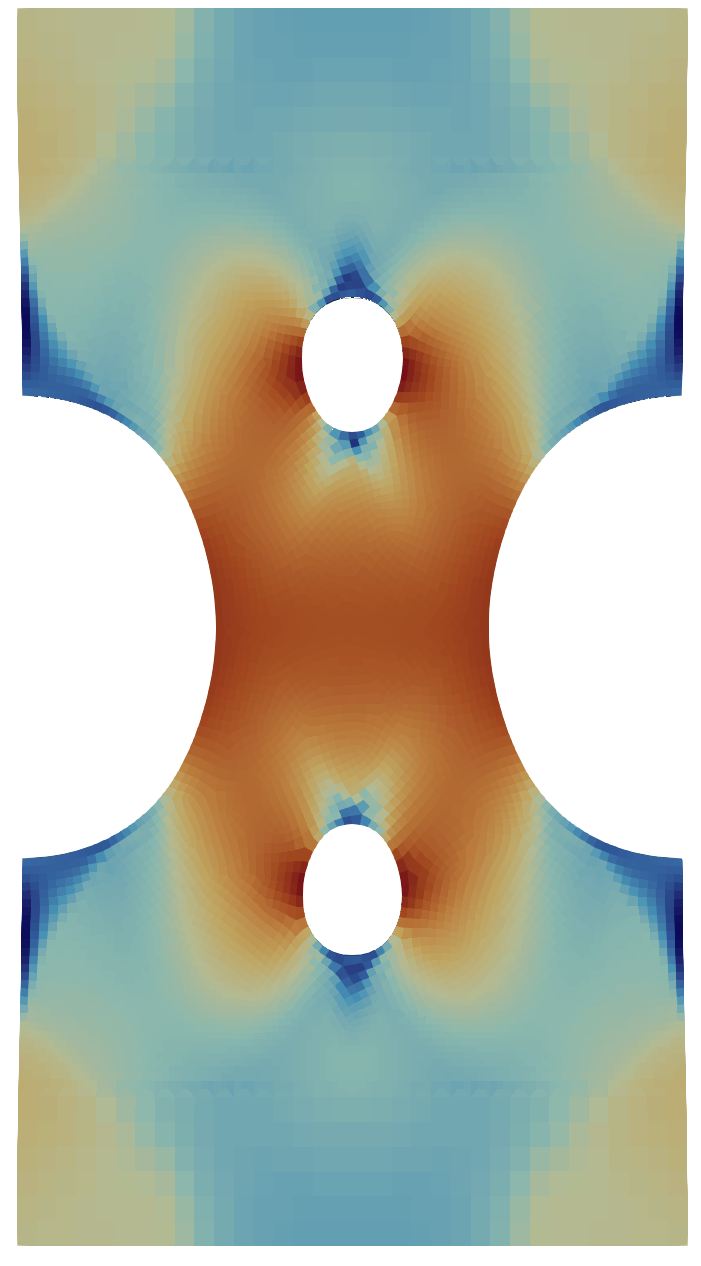}
\caption{Aluminum $\times$ specimen: (left) undeformed, (right) deformed and colored by stress along the loading direction}
\label{fig:X_schematic}
\end{figure}

To test how well the method does with typical engineering nonlinearities and measurement noise, we calibrate an elastic plastic model and treat the parameters of \eref{eq:voce} as uncertain.
We use a tension experiment for calibration and a more involved test geometry, the $\times$-specimen shown in \fref{fig:X_schematic}, for validation.
We add white noise to the calibration data but otherwise the set-up is self-consistent in the sense that the same material model is used in calibration and generating the data.

We omit the rate dependence to minimize the amount of calibration data that we need to generate to fully calibrate the model.
To generate the tension calibration data and the corresponding validation data we use a multivariate Gaussian parameter distribution:
\begin{equation}
\parameters \equiv
\{ Y_0, A, b \} \sim
\Nc(\mub, \Sigmab ) \quad \text{where} \quad
\mub = \begin{bmatrix}
34.0\\ 160.0\\ 1.95 \\
\end{bmatrix}
\quad \text{and} \quad
\Sigmab = \begin{bmatrix}
10.0&  2.0& 0.0 \\
2.0&  15.0& 0.5 \\
0.0&   0.5& 0.1 \\
\end{bmatrix}
\end{equation}
and sample sizes $n_r = 32$ for both datasets.%
\footnote{In general, we would want to compare normalized/dimensionless parameters.}
To construct the calibration data,  we added noise with a standard deviation 0.5 MPa to the data generated by the truth model.

\begin{figure}
\centering
\includegraphics[width=0.5\linewidth]{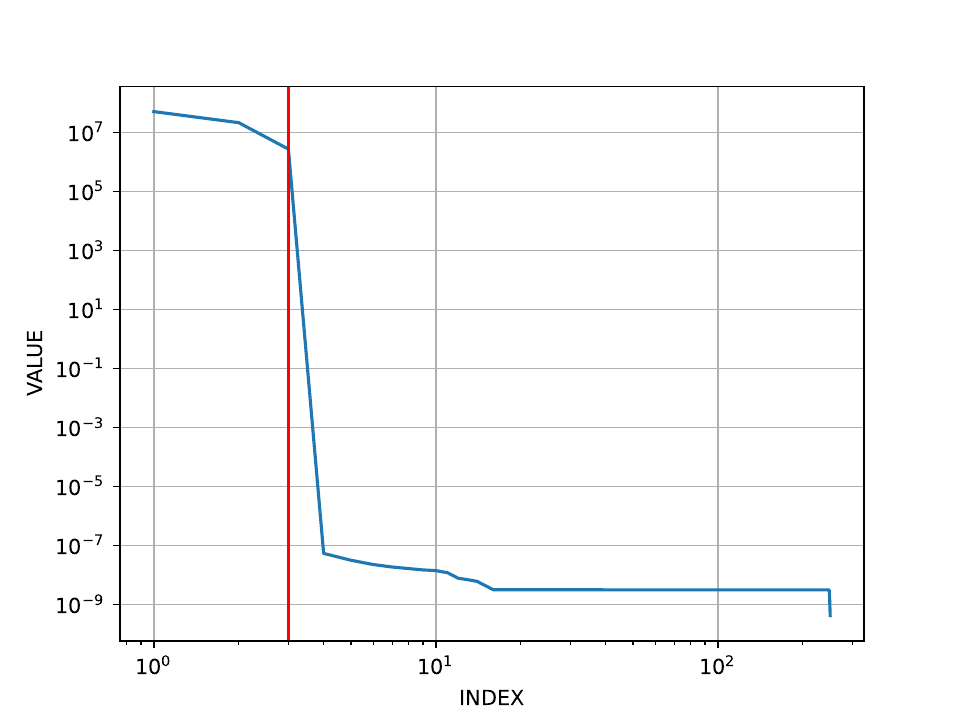}
\caption{Synthetic validation residual spectrum, where the red line indicates the number of parameters}
\label{fig:synval-singular_values}
\end{figure}

\fref{fig:synval-singular_values} shows that  residual covariance is of rank $n_\parameters$ and the singular values fall off steeply after that, which is consistent with our assumptions.

\begin{figure}
\centering
\begin{subfigure}[b]{0.57\linewidth}
\includegraphics[width=0.95\linewidth]{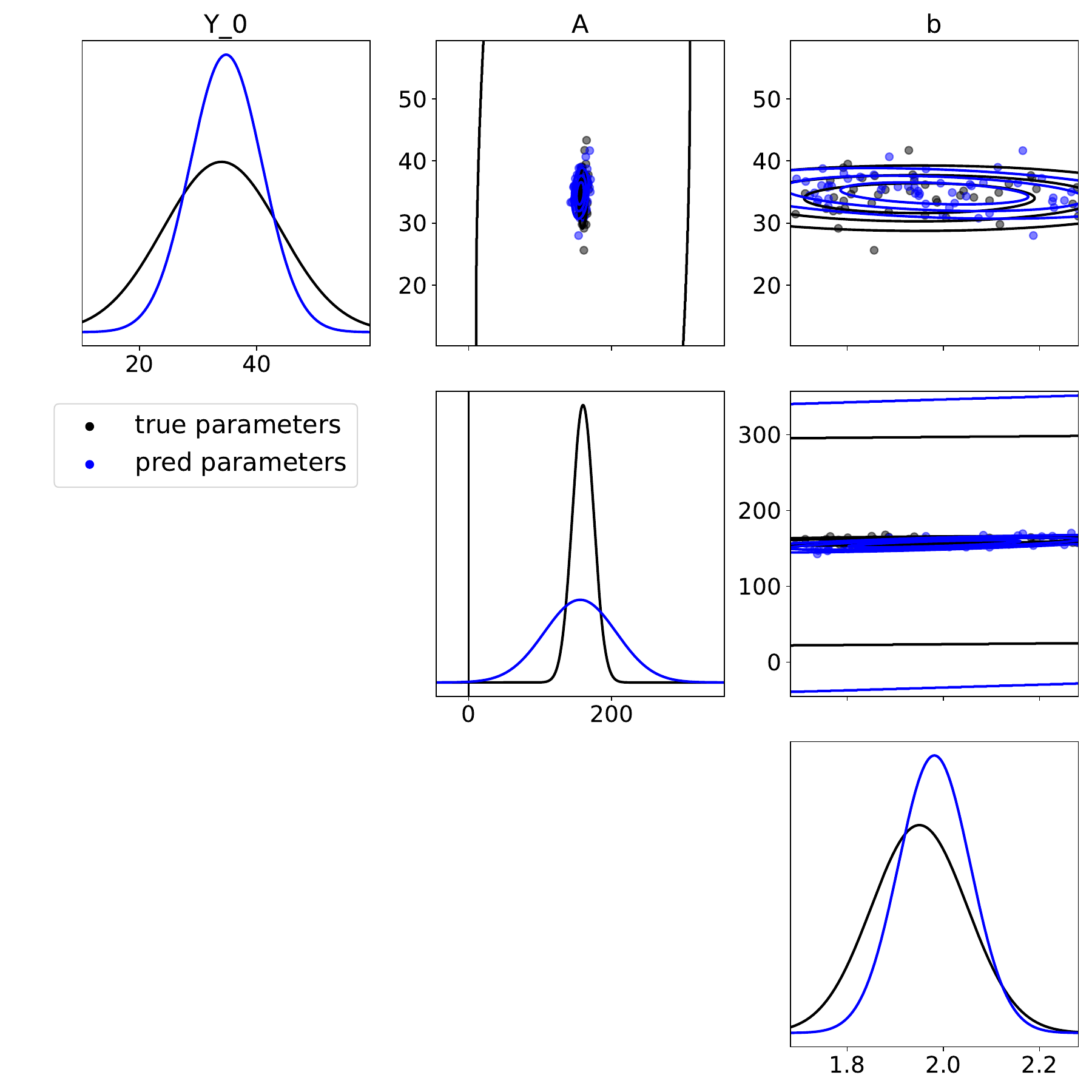}
\caption{covariance}
\end{subfigure}
\begin{subfigure}[b]{0.40\linewidth}
\includegraphics[width=0.90\linewidth]{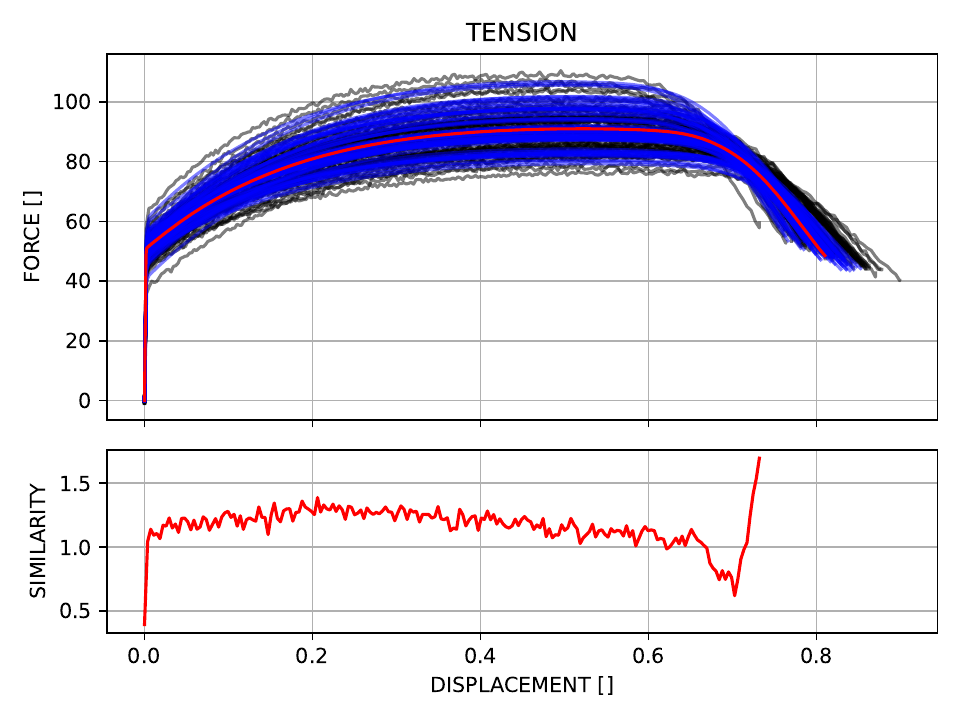}
\includegraphics[width=0.90\linewidth]{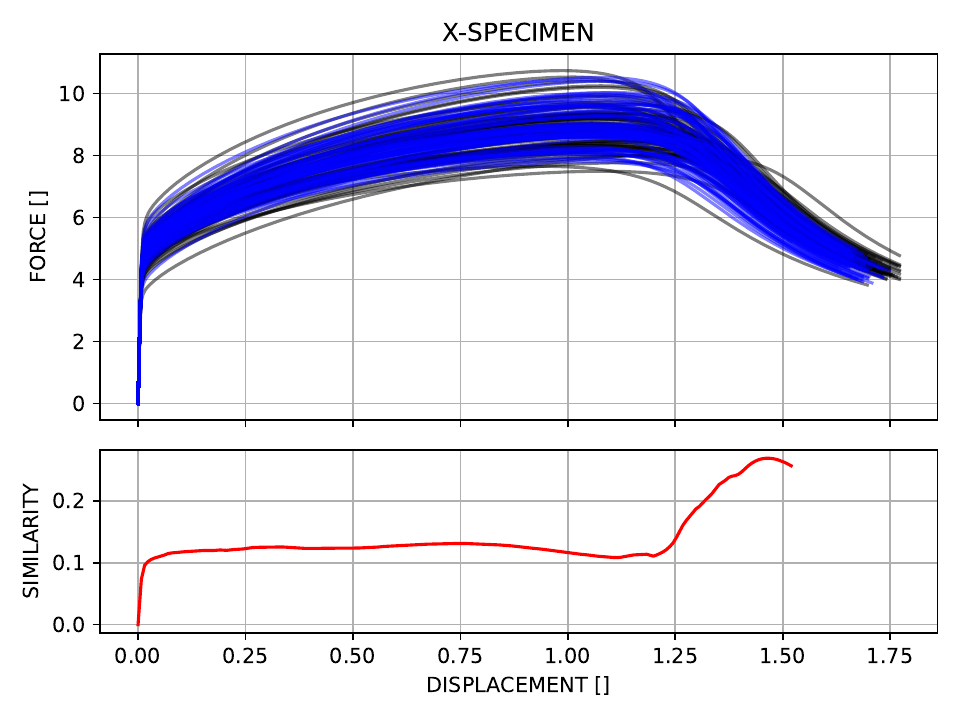}
\caption{fits}
\end{subfigure}
\caption{Synthetic validation: (a) parameter covariance, (b) push-forward predictions.
}
\label{fig:synval-fits}
\end{figure}

\fref{fig:synval-fits} shows that the true and estimated parameter distributions resemble each other but are visibly different.
Note the true and estimated means, particularly for $b$, do not coincide due to imperfect least squares calibration.
We conjecture this contributes significantly to the discrepancies in the parameter distributions.
However, the push-forwards for both the calibration and training data cover the data well.
In fact, $\Wc_1$ distance shows the similarity of the predicted and true push-forward distributions is nearly constant up until the model is pushed into extrapolation.

\subsubsection{Aluminum}
\begin{figure}
\centering
\includegraphics[width=0.45\linewidth]{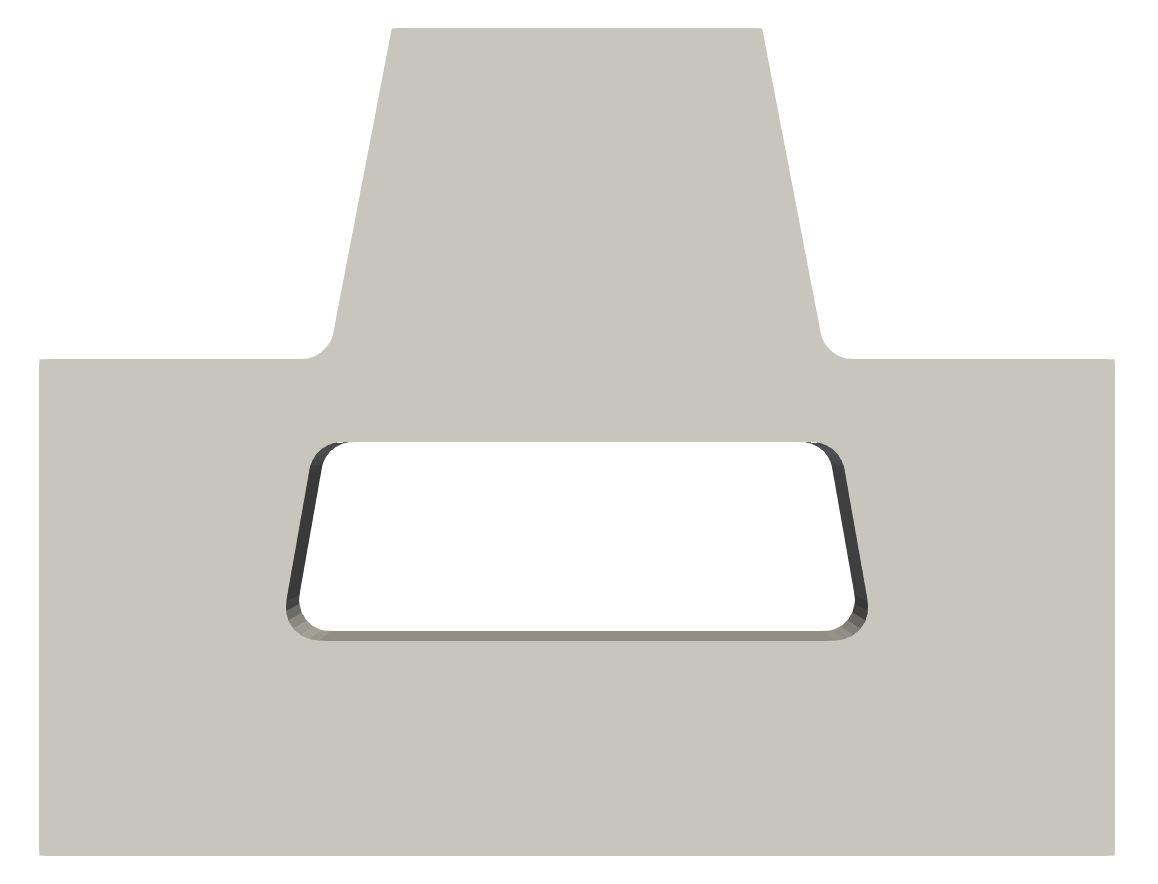}
\includegraphics[width=0.45\linewidth]{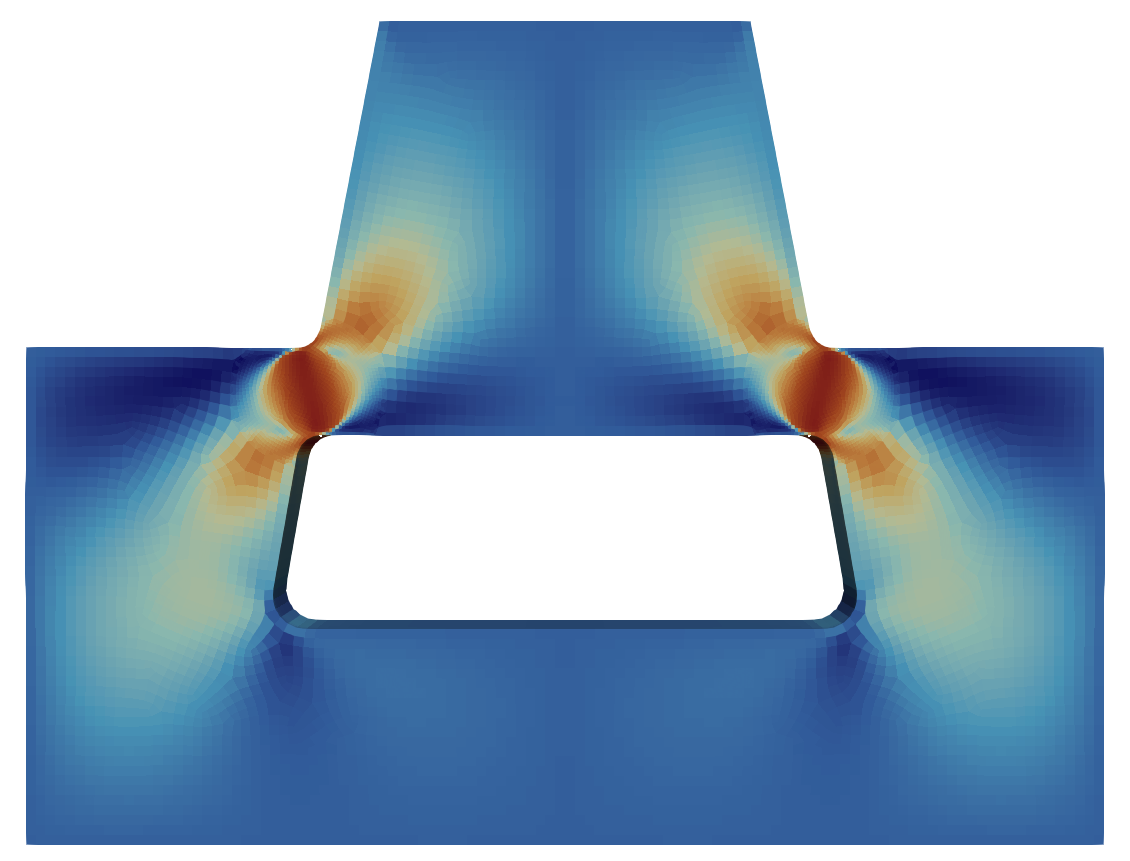}
\caption{Aluminum ``top hat'' specimen: (left) undeformed, (right) deformed and colored by shear stress}
\label{fig:aluminum_schematic}
\end{figure}

In this example, we calibrate a Hill-48 orthotropic yield surface \cite{hill1948theory} with Voce hardening \cite{voce1948relationship} to the aluminum data for a rolled plate.
One complication is that the rolling process imparts a texture to the material microstructure that results in orthotropic plasticity behavior; the experiments were designed to characterize this orthotropic behavior.
ASTM E8 uniaxial tension testing \cite{astme8} in three directions relative to the material rolling direction and a shear \emph{top-hat} testing  in six orientations relative to the material rolling direction \cite{corona2015experimental,corona2022validation}.
\fref{fig:aluminum_schematic} shows the stress and deformation of the top-hat specimen.

\begin{figure}
\centering
\includegraphics[width=0.5\linewidth]{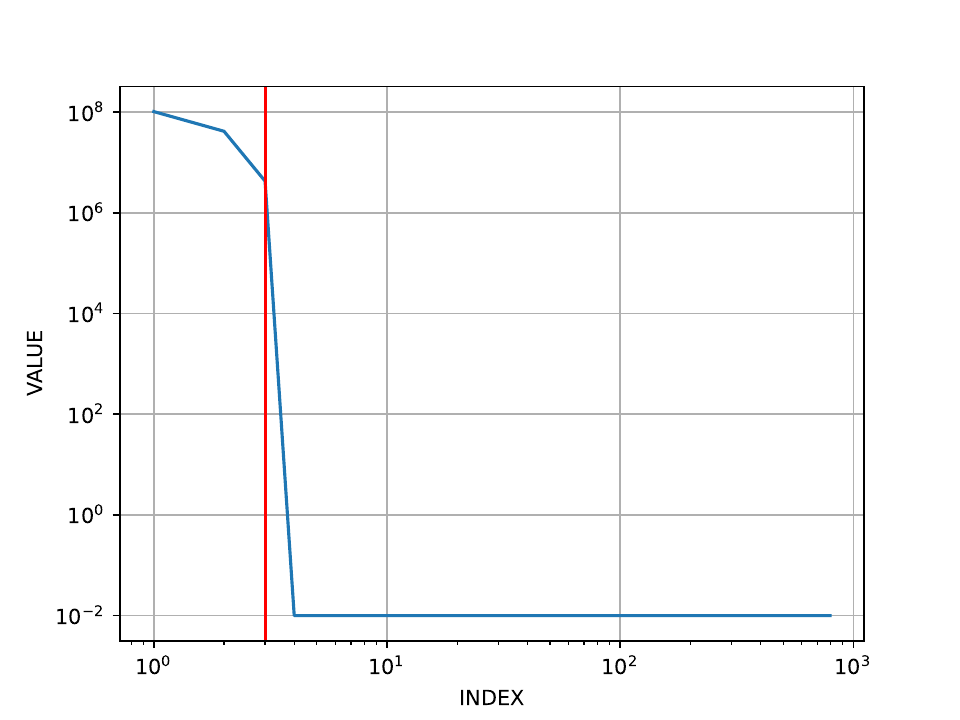}
\caption{Aluminum: residual spectrum}
\label{fig:Al_singular_values}
\end{figure}

\fref{fig:Al_singular_values} shows that  residual covariance is of rank $n_\parameters$ and the remaining singular values are consistent with the measurement noise level for these experiments.

\begin{figure}
\centering
\begin{subfigure}[b]{0.58\linewidth}
\includegraphics[width=0.95\linewidth]{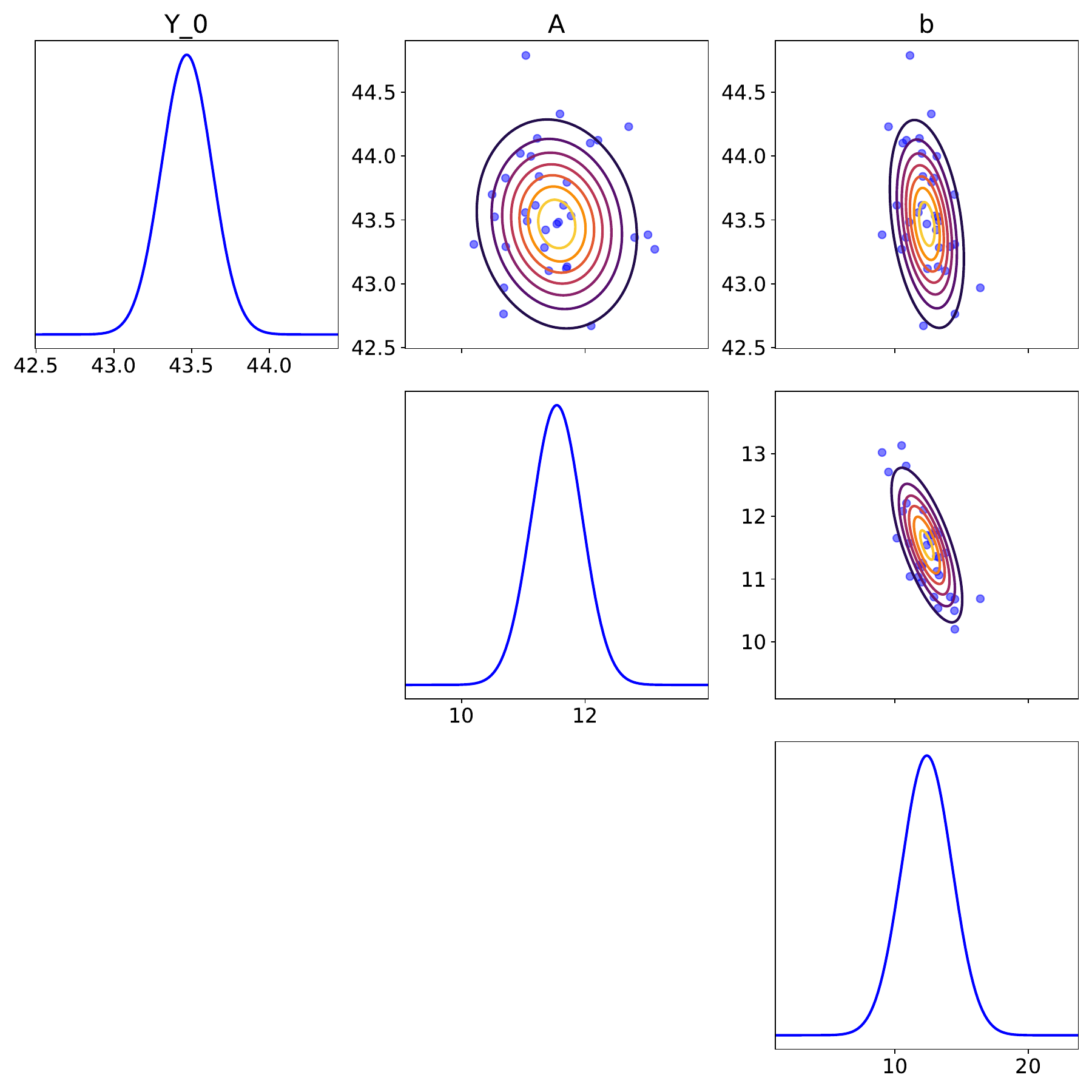}
\caption{covariance}
\end{subfigure}
\begin{subfigure}[b]{0.41\linewidth}
\includegraphics[width=0.90\linewidth]{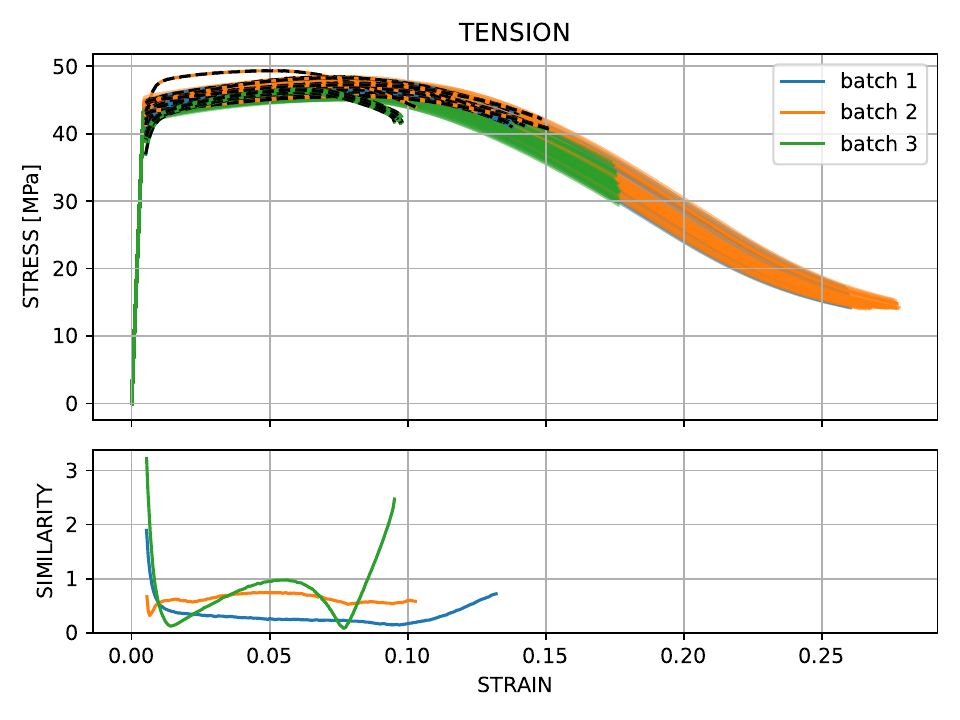}
\includegraphics[width=0.90\linewidth]{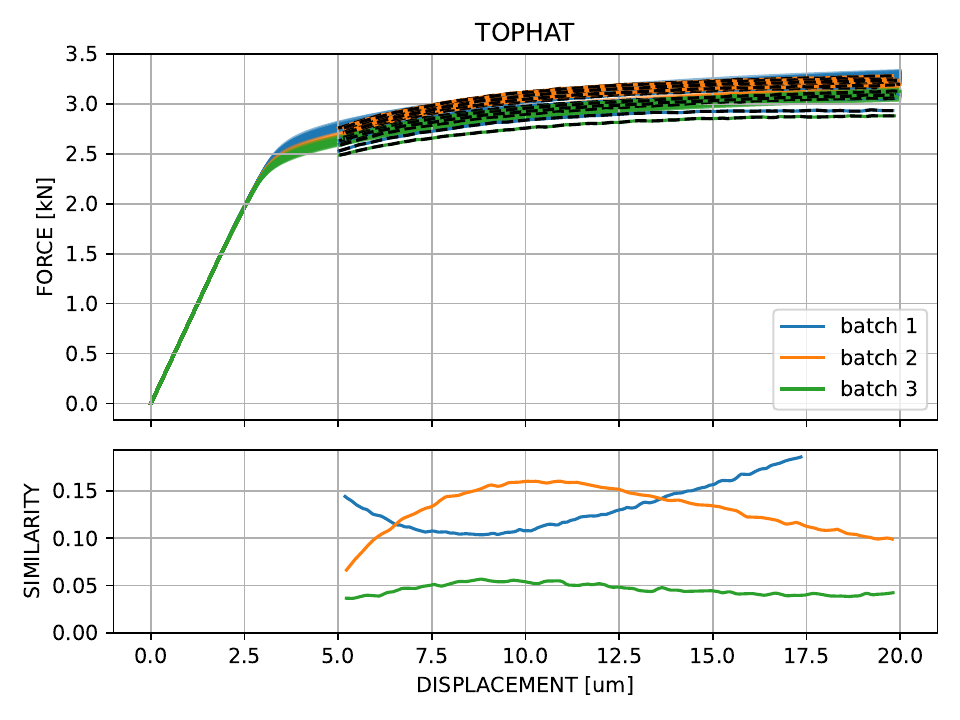}
\caption{fits}
\end{subfigure}
\caption{Aluminum calibration: (a) parameter covariance, (b) push-forward predictions.
}
\label{fig:aluminum_fits}
\end{figure}

\fref{fig:aluminum_fits}a shows that there is some correlation between $A$ and $b$.
This anticorrelation is plausible given the hardening function \eref{eq:voce} and the fact that this function largely resembles and determines the shape of the tensile response curves.
The data and predicted trajectories are shown in \fref{fig:aluminum_fits}b.
The predicted trajectories largely cover the data, with the exception of some outliers.
The $\Wc_1$ distances also indicate that the data and push-forward distributions are generally close with some exceptions near the beginning and end of the batch 1 tension data in particular.

\subsubsection{Stainless steel}
\begin{figure}
\centering
\includegraphics[width=0.65\linewidth]{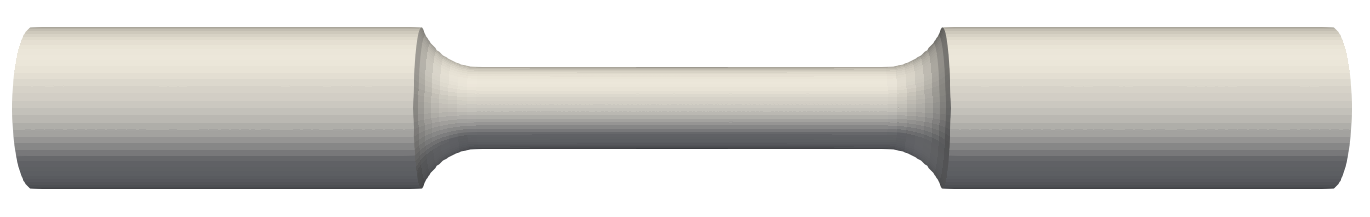}

\includegraphics[width=0.79\linewidth]{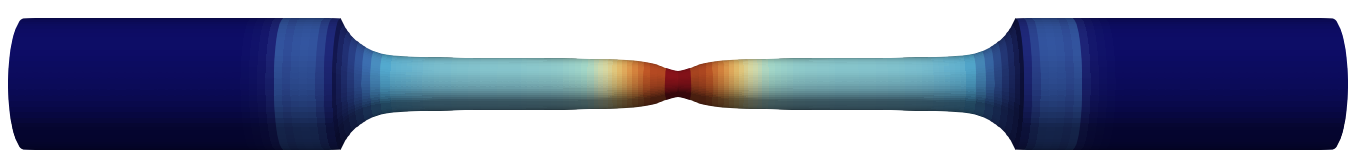}
\caption{Steel tension specimen: (upper) undeformed, (lower) deformed and colored by axial stress}
\label{fig:steel_schematic}
\end{figure}

In this example, we calibrate the rate-dependent plasticity model with a von Mises yield surface with hardening \eqref{eq:voce} to ductile failure data for 304L bar material \cite{karlson2023toward}, including  $Y_0$, $A$, $b$, and $c$.
The experimental data set used for this calibration was designed to characterize the rate-dependence of viscoplastic behavior.
A tensile bar was used for this dataset.
\fref{fig:steel_schematic} shows the undeformed shape of the standard ASTM E8 \cite{astme8} tension specimen that is held in grips at both ends and loaded in tension.
The necking/localized narrowing inducing geometric softening and the resulting axial stress field are also shown in \fref{fig:steel_schematic}.
The calibration is complicated by the variance in the failure of the specimens as indicated by black dots at the end of stress-strain trajectories.

\begin{figure}
\centering
\begin{subfigure}[b]{0.45\linewidth}
\centering
\includegraphics[width=0.95\linewidth]{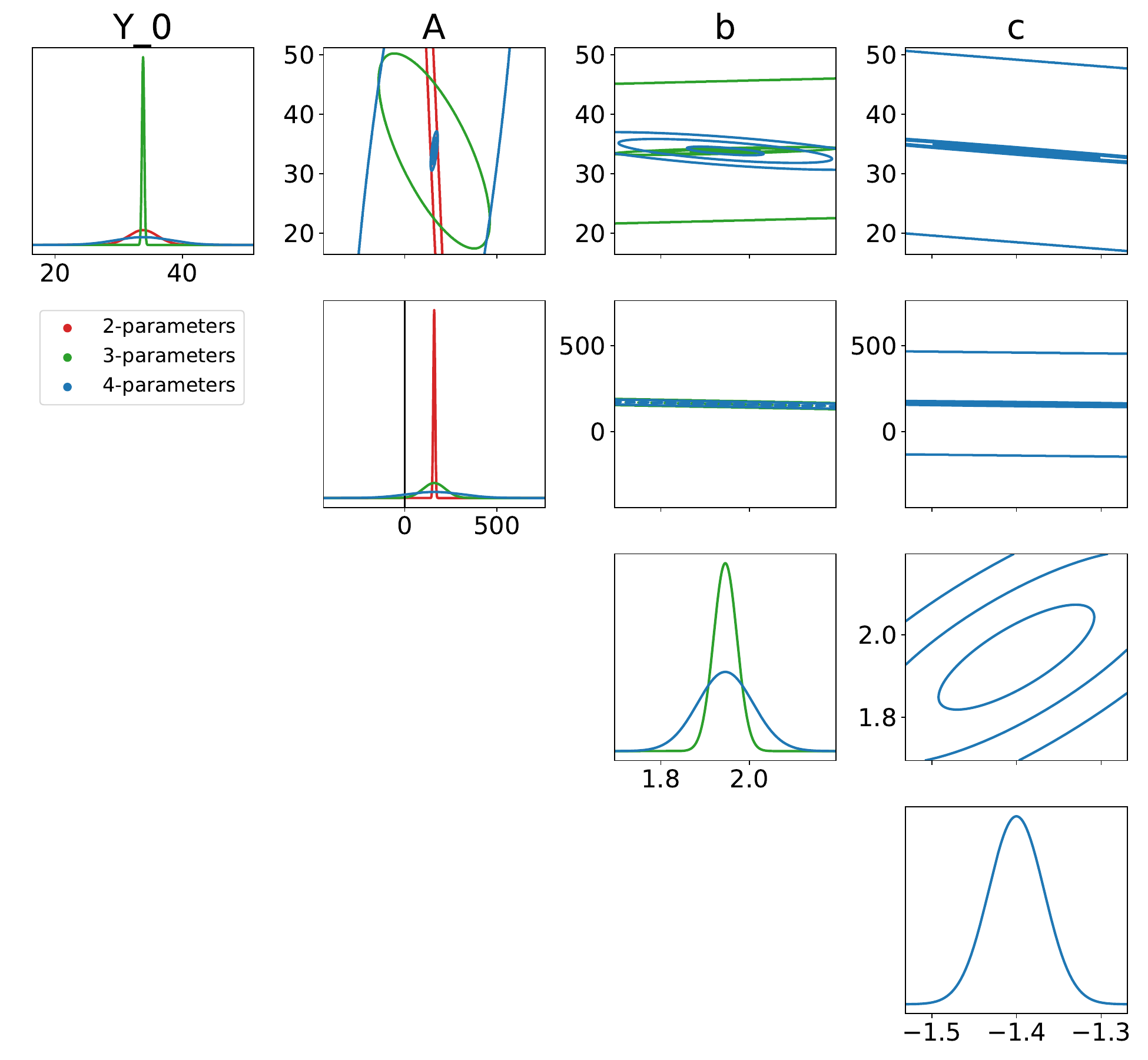}
\caption{covariance}
\end{subfigure}
\begin{subfigure}[b]{0.54\linewidth}
\centering
\includegraphics[width=0.95\linewidth]{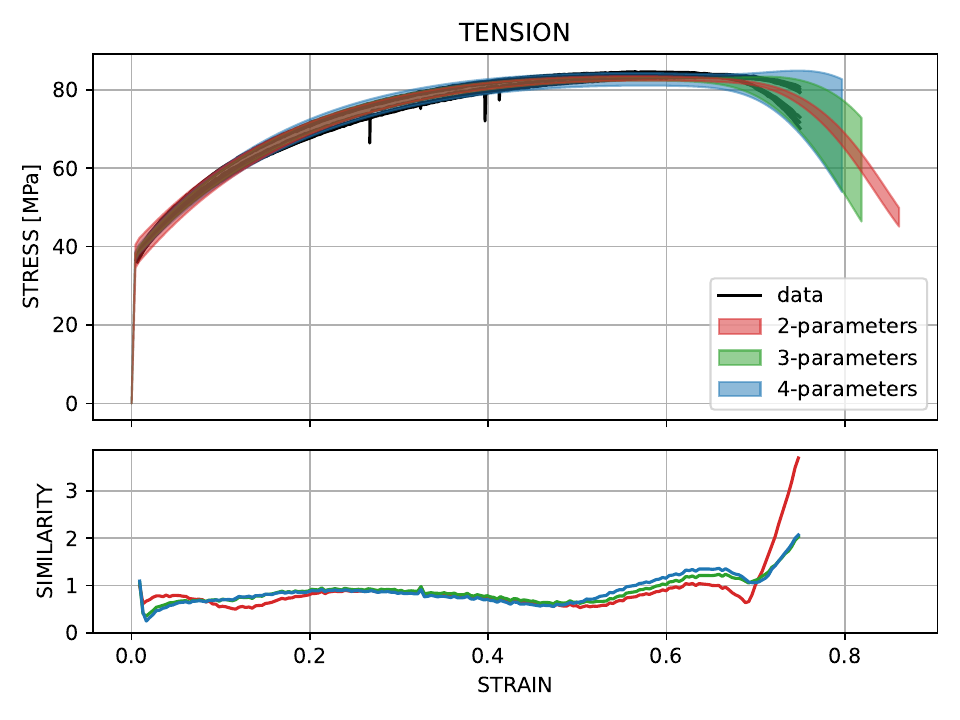}
\caption{predictions}
\end{subfigure}
\caption{Steel calibration: (a) parameter covariance, (b) push-forward predictions. Note contours on a log-scale}
\label{fig:steel_fits}
\end{figure}

We varied the number of parameters that were treated as uncertain.
\fref{fig:steel_fits}a shows a mixture of well-characterized and correlated parameters, and considerable redistribution as the size of the parameter space is changed.
Nevertheless, the push-forward trajectories, shown in \fref{fig:steel_fits}b,  clearly cover the experimental data, which implies that the systematic variations induced by the parameters are somewhat redundant in representing the tension data.
The $\Wc_1$ distances corroborate this observation.

\section{Conclusion} \label{sec:conclusion}

We proposed an extension to the classic Laplace approximation that handles both uncorrelated noise and correlated discrepancies in the mode of the \emph{embedded error} formulation \cite{sargsyan2015statistical,sargsyan2019embedded}.
This is particularly useful in calibrating material models where the material properties vary due to manufacturing inconsistencies and the like.
Furthermore, we provided an estimate of the posterior distribution using straightforward linear algebra that is often as accurate as optimization.
The method is both efficient, requiring relatively little computation after finding the MLE, and interpretable, providing a parameter covariance matrix and a noise estimate.

We plan to rectify some of the weaknesses of the method in future work.
First, the reliance on repeats to assess variance is costly in terms of experiments.
In fact, the repeated experiments currently allow us to separate correlated variance versus the uncorrelated noise.
Second, it seems that the method is fairly robust to moderate model discrepancy, but this must have limits.
We expect initializing the mean of the distribution at the MLE and then allowing it to adjust through optimization with the covariance should ameliorate some of the potential issues.
Lastly, the reliance on a multivariate normal surrogate posterior leads to violation of physical constraints, such as positivity, which is a common problem; we will explore other remedies than simple truncation and reassignment of probability mass.

\section*{Acknowledgements}
We thank Cosmin Safta (Sandia) for helpful discussions on the implementation of the BBVI method.

Sandia National Laboratories is a multi-mission laboratory managed and operated by National Technology \& Engineering Solutions of Sandia, LLC (NTESS), a wholly owned subsidiary of Honeywell International Inc., for the U.S. Department of Energy’s National Nuclear Security Administration (DOE/NNSA) under contract DE-NA0003525. This written work is authored by an employee of NTESS. The employee, not NTESS, owns the right, title and interest in and to the written work and is responsible for its contents. Any subjective views or opinions that might be expressed in the written work do not necessarily represent the views of the U.S. Government. The publisher acknowledges that the U.S. Government retains a non-exclusive, paid-up, irrevocable, world-wide license to publish or reproduce the published form of this written work or allow others to do so, for U.S. Government purposes. The DOE will provide public access to results of federally sponsored research in accordance with the DOE Public Access Plan.

\section*{Code and data availability}
The code and examples related to these developments are available at the \textsc{MatCal} GitHub site, \texttt{https://github.com/sandialabs/matcal}.

\bibliographystyle{unsrt}


\end{document}